\documentclass[11pt,notitlepage,twoside,a4paper]{article}
\usepackage{graphicx}
\usepackage{psfrag}
\usepackage{amssymb}
\usepackage{amsmath,amsthm}
\usepackage{amsmath, amssymb, amsfonts,enumerate}
\usepackage[english]{babel}
\usepackage{t1enc}
\usepackage{mathrsfs}
\newcommand {\hl }{\\\hline}
\newcommand {\ba }{\begin{array}}
\newcommand {\ea }{\end{array}}
\newcommand{\ds}{\displaystyle}
\newcommand{\dist}{\mbox{\,\,\rm dist}}

\newtheorem{defi}{Definition}[section]
\newtheorem{theo}[defi]{Theorem}
\newtheorem{propo}[defi]{Proposition}

\newtheorem{lem}[defi]{Lemma}

\newtheorem{rem}[defi]{Remark}

\newcommand{\labell}[1]{\label{#1}}
\author{\\ \small IRMA - UMR 7501 CNRS/ULP\\\small 7 rue René
  Descartes - 67084 Strasbourg Cedex, France\\\small email: 
sellama@math.u-strasbg.fr}
\title{\bf On the distance between separatrices for the  discretized pendulum equation  }
\date{}
\numberwithin{equation}{section}
\begin{document}
\maketitle
 \renewcommand{\abstractname}{abstract}
  \begin{abstract}
\noindent We consider the discretization
\begin{equation*}
q(t+\varepsilon)+q(t-\varepsilon)-2q(t)=\varepsilon^{2}\sin\big(q(t)\big),
\end{equation*}
%This second order equation is a discretization of the  pendulum equation  $ q '' = \sin (q) $. It is equivalent to  the following system of  first order difference equations
%y(t+\varepsilon)=y(t-\varepsilon)+2\varepsilon\big(1-y(t)^{2}\big),
%\end{equation*}
$\varepsilon>0$ a small parameter, of the pendulum equation  $ q '' = \sin (q) $;
in system form, we have the discretization
\begin{equation*}
q(t+\varepsilon)-q(t)=\varepsilon p(t+\varepsilon),
  \ p(t+\varepsilon)-p(t)=\varepsilon\sin\big(q(t)\big).
\end{equation*}
of the system
\begin{equation*}
q'=p,\ p'=\sin(q).
\end{equation*}

The latter system of ordinary differential equations
has two saddle points at  $A=(0,0)$, $B=(2\pi, 0)$ and
near both, there exist stable and unstable manifolds. It also admits a 
heteroclinic orbit connecting the stationary points $B$ and $A$ parametrised by
$q_0(t)=4\arctan\big(e^{-t}\big)$ and which contains the stable manifold of this system at $A$ as well as its unstable manifold at $B$. 
We prove that the stable
manifold  of the point $A$ and the unstable
manifold of the point $B$ do not coincide for the discretization. 
More precisely, we show that the vertical distance  between these two manifolds is exponentially small but not zero and in particular we give an asymptotic estimate of this distance.  For this purpose we use a method adapted from the article of Sch\"afke-Volkmer \cite{SV} using formal series  and  accurate  estimates  of the coefficients. Our result is similar to that of  Lazutkin et.\ al.\ \cite{LS}; our method of proof, however,  is quite different.

 % In this paper we will use a method adapted from the  paper of
 % Sch\"afke-Volkmer \cite{SV} using formal series  and  accurate  estimates  of the coefficients to show that  the stab%le manifold $W_{s}^{+}$\, of the point $A=(1,1)$ and the
  % unstable manifold  $W_{i}^{-}$\, of the point $B=(-1, -1)$ for the discretized   logistic equation
 % do not coincide moreover the distance in the sense of Hausdorff between these two manifolds is exponentially small an%d not zero, in particular we give an
% asymptotic approximate of this distance.
 \end{abstract}

\section{Introduction}
We consider the following difference equation
\begin{eqnarray}
q(t+\varepsilon)+q(t-\varepsilon)-2q(t)=\varepsilon^{2}\sin\big(q(t)\big).
%\begin{cases}
%q(t+\varepsilon)+q(t-\varepsilon)-2q(t)=\varepsilon^{2}\sin\big(q(t)\big),\\
 %q(t)\to 0 \,\,\,\, \text{As} \,\,\,\, t\to -\infty,\ q(0)=\pi.
%\end{cases}
\labell{diff}\end{eqnarray}
This second order equation is a discretization of the   pendulum equation  $ q '' = \sin (q) $. It is equivalent to  the following system of  first order difference equations
\begin{eqnarray}
\begin{cases}
q(t+\varepsilon)=q(t)+\varepsilon p(t+\varepsilon),\\

 p(t+\varepsilon)=p(t)+\varepsilon\sin\big(q(t)\big).
\end{cases}
\labell{int1}\end{eqnarray}
which can be considered as a discretization of the system
\begin{eqnarray}
\begin{cases}
q'=p,\\
p'=\sin(q).
\end{cases}
\labell{int3}\end{eqnarray}
The latter system has two saddle points at  $A=(0,0)$, $B=(2\pi, 0)$ and
there exist stable and unstable manifolds.  For the discretized equation (\ref{int1}) and sufficiently small $\varepsilon>0$, these manifolds still exist.  

The system (\ref{int3}) has $\big(q_0(t),q^{'}_0(t)\big)$, where $q_0(t)=4\arctan\big(e^{-t}\big)$, as a heteroclinic orbit connecting the stationary points $B$ and $A$; it is a parametrisation of  the curve $p=-2\sin(q/2)$ and contains the stable manifold of (\ref{int3}) at the point $A$ as well as its unstable manifold at $B$. This curve, together with $p=2\sin(q/2)$, separates regions with periodic orbits from regions with non-periodic orbits and is therefore often called a {\em separatrix}.
Our purpose is study the behavior of this separatrix under discretization of the equation -- it turns out that there is no longer a heteroclinic orbit for system (\ref{int1}) and its the stable manifold at $A$ and the unstable manifold at $B$ no longer coincide. More precisely, we want to estimate the distance  between the stable manifold $W_{s, \varepsilon}^{-}$  of (\ref{int1}) at $A$ and the unstable manifold    $W_{u, \varepsilon}^{+}$  of (\ref{int1}) at $B$ as a function of the parameter $\varepsilon$.

Lazutkin et.\ al.\ \cite{LS}, Gelfreich \cite{G}, (see also Lazutkin \cite{L1}\cite{L2}) had given an asymptotic estimate of the splitting angle between the  manifolds. Starting from a heteroclinic solution
% $q_{1}(t)=4\arctan\big(e^t\big)$ 
of the differential equation, they study the behavior of analytic solutions of the  difference equation  in the neighbourhood of its singularities $t=\pm\frac{\pi}{2}i$.
%using a method which primarily consists to study the analytical solutions of the  difference equation in the neighbourhood of singularities of the solution of the differential equation.

We show that the distance 
between these two manifolds is exponentially small but not zero and we give an asymptotic estimate of this distance. This result is similar to that of  Lazutkin et.\ al.\ \cite{LS}; our method of proof, however,  is quite different.

We use a method adapted from the article of Sch\"afke-Volkmer \cite{SV} using a formal power series solution and accurate estimates of the coefficients. This method was adapted for the logistic equation in Sellama\cite{S}.  It turns out that the adaptation of this method for the  pendulum equation  is more difficult than in the case of the logistic equation. 

%Fruchard-Sch\"afke \cite{FS}(subsection 5.2) had shown that there are two families of entire functions $Y_{\varepsilon}^{\pm}(t)=(q^{\pm}_{\varepsilon}(t),p^{\pm}_{\varepsilon}(t))$ solutions of  system
%(\ref{int1}) such that the functions $t\mapsto\big(q^{+}_{\varepsilon}(t),p^{+}_{\varepsilon}(t)\big) $, 
% $t\mapsto\big(q^{-}_{\varepsilon}(t),p^{-}_{\varepsilon}(t)\big) $ provide parametrizations $w_{u, \varepsilon}^+$ of  $W_{u, \varepsilon}^+$ respectively $w_{s, \varepsilon}^-$ of  $W_{s, \varepsilon}^-$. As the system (\ref{int1}) is autonomous, replacing $t$ by
%$t+\gamma$ with constant $\gamma$ provides another parametrisation. We choose it such 
%that $q^-_\varepsilon(0)=\pi$ respectively $q^+_\varepsilon(0)=\pi$.

We will show 
 %Indeed, to find an estimate asymptotic of the coefficients of the formal solution we were confronted with the difficulty in twice reversing the operator on the polynomials. 
\begin{theo} Given any positive $t_0$, it is known that for sufficiently smal  $\varepsilon_0>0$ and all $t\in ]-\infty, t_0] $ there is exactly one one point $w^{+}_{u,\varepsilon}(t)=(q_0(t),\tilde p^+_{u,\varepsilon}(t))$  on the stable unstable manifold having first coordinate $q_0(t)$.
There exist constants $\alpha\neq0$, such that for any positive $t_0$ 
\begin{equation*} 
\dist_{v}\big(w^{+}_{s,\varepsilon}(t),W^{-}_{s,\varepsilon}\big)=\dfrac{4\pi\alpha}{\varepsilon^{2}}\cosh(t)\sin\Big(\frac{2\pi t}{\varepsilon}\Big)e^{-\frac{\pi^{2}}{\varepsilon}}+O\bigg(\frac{1}{\varepsilon}e^{-\frac{\pi^{2}}{\varepsilon}}\bigg),\quad \text{as} \,\,\,\varepsilon\searrow 0,
\end{equation*}
uniformly for $-t_0<t<t_0$, 
where$\dist_{v}\big(P,W^{-}_{u,\varepsilon}\big)$ denotes the vertical distance of a point $P$ from the unstable manifolds $W^{-}_{u,\varepsilon}$.
\end{theo}
This result corresponds to the result of Lazutkin et.\ al.\ \cite{LS} 
as the angle between the manifolds  at an intersection point
is asymptotically equivalent to  
 
 \noindent$\tfrac1{q_0'(t)}\tfrac d{dt}\dist_{v}\big(w^{+}_{s,\varepsilon}(t),W^{-}_{u,\varepsilon}\big)$, but we do not want to give any detail here.
% where $q_0(t)=4\arctan\big(e^t\big)$ is  a heteroclinic solution of the differential equation

Our proof uses the following steps. First, we construct a  formal solution for the difference equation  (\ref{diff}) in the form of a power series in $d=2\rm{\mbox{arsinh}}(\varepsilon/2)$, whose  coefficients are polynomials in $u=\tanh(d t/\varepsilon)$. This is done in section 2; the introduction of $d$ is necessary because  polynomials are desired as coefficients. Then, we
 give  asymptotic approximations  of these coefficients using appropriate norms on  spaces of polynomials.  To that purpose we introduce operators on  polynomials series. In section 6 we use the truncated Laplace transform to construct a function which satisfies (\ref{diff}) except for an   exponentially small error. The next and last step is to give an asymptotic estimate for the distance of some point of   the stable manifold from the unstable manifold. 
A calculation shows that  $\alpha=89.0334$ and therefore $4\pi\alpha=1118.8267$ (See Remark \ref{rem1}); the corresponding constants of Lazutkin have already been calculated with high precision (See Lazutkin et.\ al.\ \cite{LS}).
A {\em proof} that $\alpha\neq0$ as in \cite{SV} or \cite{S} would be possible.  Y.B.  Suris \cite{SY} had shown that $\alpha\neq 0$.
%but seems not te be too interesting.

%The next step is to find out the asymptotic behaviorof the sequence of these cofficients.

\section{Formal solutions}
%Let us  consider the following  difference equation
%{ \bf DIFF}
%\begin{eqnarray}
%\begin{cases}
%q(t+\varepsilon)+q(t-\varepsilon)-2q(t)=\varepsilon^{2}\sin\big(q(t)\big),\\
% q(t)\to 0 \,\,\,\, \text{as} \,\,\,\, t\to -\infty,\ q(0)=\pi.
%\end{cases}
%\labell{diff}\end{eqnarray}
 The purpose of this section is to find a convenient formal solution  for equation (\ref{diff}). First, we need some preparations. We put
\begin{eqnarray*}
u:&=&\tanh\bigg(\frac{d}{\varepsilon}t\bigg),\\
q_{0d}(t):&=&4\arctan\Bigg(\exp\bigg(-\frac{d}{\varepsilon}t\bigg)\Bigg),\\
q_d(t)&=&\sqrt{1-u^{2}}A_{d}(u)+q_{0d}(t), \ \ A_d(u)=\sum_{n=1}^\infty A_n(u)d^n
\end{eqnarray*}
for a formal solution of (\ref{diff}),
where $d=\varepsilon+\sum_{n=3}^{\infty}d_{n}\varepsilon^{n}$ is a formal powers series in $\varepsilon$ to be determined.

{\bf Remark.} The  linearization of  equation (\ref{diff}) at the point $A$    gives the following equation
\begin{eqnarray*}
Z(t+\varepsilon)+Z(t-\varepsilon)-2Z(t)=\varepsilon^{2}Z(t).
\end{eqnarray*}
The parameter $d$ is such that $Z(t)=e^{-dt}$ is a solution of this equation,  therefore $\varepsilon$ and $d$ are coupled by the relation $d=2\text{arcsinh}(\varepsilon /2).$

By  Taylor expansion, we obtain  
\begin{eqnarray}
q_{0d}(t+\varepsilon)+q_{0d}(t-\varepsilon)-2q_{0d}(t)=2\sum^{+\infty}_{n=1}\frac{1}{(2n)\,!}q^{(2n)}_{0d}(t) \,\varepsilon^{2n},
\end{eqnarray}
where $\dfrac{2}{(2n)!}q^{(2n)}_{0d}(t)\varepsilon^{2n}/d^{2n}$ is an odd polynomial 
$I_{2n-1}(u)$  multiplied by $\sqrt{1-u^{2}} $; we find $I_{2n-1}(1)=4/(2n)!$.
%, since  $\dfrac{q^{(2n)}_0(t)}{\sqrt{1-u^{2}}}\to 0$ as $t \to \infty$, we have  $I_{2n-1}(1)=0$.

Using  $\cos(q_{0d})=2u^{2}-1,\,\,\,
\sin(q_{0d})=2u\sqrt{1-u^{2}}$,  we can express our equation (\ref{diff}) in the form 

\begin{eqnarray}
%\begin{cases}
A_{d}(T^{+})\sqrt{\frac{1-(T^{+})^{2}}{1-u^{2}}}+A_{d}(T^{-})\sqrt{\frac{1-(T^{-})^{2}}{1-u^{2}}}-2A_{d}(u)=f\big(\varepsilon,u,A_{d}(u)\big)%\\ A_{d}(0)=0
%\end{cases}
\labell{eq1}\end{eqnarray}
or equivalently
\begin{eqnarray}
%\begin{cases}
\dfrac{A_{d}(T^{+})}{\cosh(d)+u\sinh(d)}+\dfrac{A_{d}(T^{-})}{\cosh(d)-u\sinh(d)}-2A_{d}(u)=f\big(\varepsilon,u,A_{d}(u)\big)%\\ A_{d}(0)=0
%\end{cases}
\labell{eq2}
\end{eqnarray}
where 
\begin{eqnarray*}
f\big(\varepsilon,u,A_{d}(u)\big)&=&\varepsilon^{2}\bigg(2u\cos\Big(A_{d}(u)\sqrt{1-u^{2}}\Big)+\frac{2u^{2}-1}{\sqrt{1-u^{2}}}\sin\Big(A_{d}(u)\sqrt{1-u^{2}}\Big)\bigg)\\
&-&
\sum_{n=1}^{+\infty}I_{2n-1}(u)\, d^{2n},\\
T^{+}&=&T^{+}(d,u)=\frac{u+\tanh(d)}{1+u\tanh(d)}=\tanh\Big(\dfrac{d}{\varepsilon}(t+\varepsilon)\Big),\\
T^{-}&=&T^{-}(d,u)=\frac{u-\tanh(d)}{1-u\tanh(d)}=\tanh\Big(\dfrac{d}{\varepsilon}(t-\varepsilon)\Big).
%d:&=&\sum^{+\infty}_{n=1}{d_{n}\varepsilon^{n}}.
\end{eqnarray*}
As $u\rightarrow1,$ the expressions $T^+$ and $T^-$ reduce to 1, the denominators in (\ref{eq2}) simplify to $e^{\pm d}$ and hence equation
(\ref{eq2}) reduces to $$(e^{-d}+e^d-2)A_d(1)=\varepsilon^2(2+A_d(1))-4(\cosh(d)-1).$$
This is equivalent to  $(2\cosh(d)-2-\varepsilon^2)(2+A_d(1))=0$ and hence we have
necessarily $\varepsilon=2\sinh(d/2)$ if we want a formal solution such that the coefficients have limits as $u\rightarrow1$.

%We are looking for  formal solutions in the form of formal power series, generally this  is done by  finding a formula of recurrence using identification of coefficients. 
\begin{theo}\textbf{(On the formal solution)} If
  $\varepsilon=2\sinh(d/2)$, 
 then equation (\ref{eq1}) has a unique formal  solution of the form 

\begin{equation}
A_{d}(u)=\sum^{+\infty}_{n=1}A_{2n-1}(u)d^{2n},
\end{equation}
where $A_{2n-1}(u)$ are odd polynomials of degree $\leq 2n-1$.
%   and $A_{2n-1}(0)=0$ for $n\geq 1$.
\labell{the1}\end{theo}
\noindent{\bf Remark:}
A similar  formal solution was found using another method in  \cite{SY}.
\noindent\textbf{Proof.}
 We will use the Induction Principle  to show that there
 exist unique odd polynomials $A_{1},A_{3},A_{5}...A_{2n-1}$ such that 
  \begin{equation}
Z_{n}(d,u)=\sum_{k=1}^{n}A_{2k-1}(u)d^{2k} 
\end{equation} 
satisfy

\begin{equation}
R_{n}(d,u)=O(d^{2n+4})
\labell{re1}\end{equation}
where
\begin{eqnarray}
R_{n}(d,u)&=&Z_{n,d}\big(T^{+}\big)\sqrt{\frac{1-(T^{+})^{2}}{1-u^{2}}}+Z_{n,d}\big(T^{-}\big)\sqrt{\frac{1-(T^{-})^{2}}{1-u^{2}}}\nonumber\\
&-&2Z_{n,d}(u)-f\big(\varepsilon,u,Z_{n,d}(u)\big)
\labell{re2.7}\end{eqnarray}

\noindent For $n=1$, a short calculation shows that we must have 
$A_{1}(u)=-\frac{1}{4}u $ and hence $Z_{1,d}(u)=-\frac{1}{4}ud^{2}$. We obtain 
\begin{equation*}
R_{1}(d,u)=(\frac{-91}{48}u^{5}+\frac{137}{48}u^{3}-\frac{23}{24})d^{6}+O(d^{8}).
\end{equation*}
 Suppose now that there exists $A_{1},A_{3},A_{5}...A_{2n-1}$ such that
\begin{equation}
Z_{n}(d,u)=\sum_{k=1}^{n}A_{2k-1}(u)d^{2k}
\end{equation}
 satisfies (\ref{re1}), (\ref{re2.7}). We show that there is
a unique polynomial $A_{2n+1}(u)$ such that 
\begin{equation}
 Z_{n+1}(d,u)=Z_{n}(d,u) +A_{2n+1}(u)d^{2n+2}
\labell{re2}\end{equation}
 satisfies (\ref{re1}).
We put 
\begin{equation}
R_{n}(d,u)=R_{2n+3}(u)d^{2n+4}+O\big(d^{2n+6}\big) 
\labell{re3}\end{equation}
where $R_{2n+3}(u)$  is odd and $\text{deg}(R_{2n+3}(u))\leq2n+3$.
%and  $R_{2n+3}(0)=0$.

We substitute $ Z_{n+1}(d,u)$ in equation (\ref{re2.7}). Using Taylor expansion, (\ref{re2}), (\ref{re3}) and  $\varepsilon=2\sinh(d/2)$, we obtain 
\begin{eqnarray*}
&&Z_{n+1,d}\big(T^{+}\big)\sqrt{\frac{1-(T^{+})^2}{1-u^{2}}}-Z_{n+1,d}\big(T^{-}\big)\sqrt{\frac{1-(T^{-})^2}{1-u^{2}}}-2Z_{n+1,d}(u)-\\
&& f\big(\varepsilon,u,Z_{n+1,d}\big)=\bigg[(u^{4}-2u^{2}+1)A^{''}_{2n+1}(u)+(4u^{3}-4u)A^{'}_{2n+1}(u)+\\
&& R_{2n+3}(u)\bigg]d^{2n+4}+ O\big(d^{2n+6}\big)
\end{eqnarray*}
We  notice that (\ref{re3}) is satisfied if only if
\begin{equation}
\big[(1-u^{2})^{2}A^{'}_{2n+1}(u)\big]^{'}+R_{2n+3}(u)=0
\end{equation}

This differential equation 
has a unique   solution vanishing at $u=0$ without singularity at $u=1$, namely 
\begin{equation}
A_{2n+1}(u)=-\int_0^{u}\frac{\int_{1}^{t}R_{2n+3}(s)ds}{(1-t^{2})^{2}}dt.
\end{equation}
We now show that this solution is an odd polynomial  of $u$. It is clear that $\int_{1}^{t}R_{2n+3}(s)ds $ vanishes for $t=1$ and as $R_{2n+3}(s)$  is odd,  it  also vanishes for $t=-1$. It suffices to show that $R_{2n+3}(s) $ also vanishes at $t=\pm 1 $. Indeed, taking  the limit of   (\ref{re2.7}) as   $u\to 1$ as we did for (\ref{eq2}) and using 
\begin{eqnarray*}
\lim_{u\to1} f\big(\varepsilon,u,Z(d,u)\big)&=&\varepsilon^2 Z(d,1)
%\\ \lim_{u\to1}\sqrt{\frac{1-(T^{+})^2}{1-u^{2}}}&=&\sqrt{\dfrac{1-\tanh(d)}{\tanh(d)+1}}\\ \lim_{u\to1}\sqrt{\frac{1-(T^{-})^2}{1-u^{2}}}&=&\sqrt{\dfrac{1+\tanh(d)}{1-\tanh(d)}}
\end{eqnarray*}
 we obtain 

\begin{equation*}
%R_{2n+3}(1)d^{2n+4}=\Bigg(\sqrt{\dfrac{1-\tanh(d)}{\tanh(d)+1}}+\sqrt{\dfrac{1+\tanh(d)}{1-\tanh(d)}}-2-\varepsilon^{2}\Bigg)Z(d,1)+O(d^{2n+6}).
R_{2n+3}(1)d^{2n+4}=\Bigg(e^d+e^{-d}-2-\varepsilon^{2}\Bigg)Z(d,1)+O(d^{2n+6}).
\end{equation*}
By our choice of $\varepsilon=2\sinh(d/2)$, we obtain 
 $R_{2n+3}(1)d^{2n+4}=O(d^{2n+6})$.\, Consequently  $R_{2n+3}(1)=0$.
As $R_{2n+3}(u)$ is odd, we also have $R_{2n+3}(-1)=-R_{2n+3}(1)=0$. This proves that 
  $A_{2n+1}(u)$ is an odd polynomial of  $\text{degree}\big(A_{2n+1}(u)\big)\leq2n+1$
 and  $A_{2n+1}(0)=0.$

The first polynomials  $A_{2n-1}(u)$ with  $n>0$  can be calculated using Maple.
$$\begin{array}{|c|c|c|c|c|} \hline

n&1 &2 &3     \hl

 A_{2n-1}(u)\  &-\frac{1}{4}u &\bigg(\frac{91}{864}u^{3}-\frac{47}{576}u \bigg)
 &\bigg(-\frac{319}{2880}u^{5}+\frac{185}{1152}u^{3}-\frac{3703}{69120}u \bigg) \hl

\end{array}
$$

Now, we ntroduce the operators $\mathcal{C}_{2}, \mathcal{C},
\mathcal{S}_{2}, \mathcal{S}$ defined by
 %In order to rewrite equation (\ref{eq2}), 
\begin{equation}          
\begin{array}{llll}
    \mathcal{C}(Z)(d,u)=\frac{1}{2}\big(Z(d,T^{+\frac{1}{2}})+Z(d,T^{-\frac{1}{2}})\big)&\\
\\
    \mathcal{S}(Z)(d,u)=\frac{1}{2}\big(Z(d,T^{+\frac{1}{2}})-Z(d,T^{-\frac{1}{2}})\big)&\\
\\
\mathcal{C}_{1}(Z)(d,u)=\frac{1}{2}\big(Z(d,T^{+})+Z(d,T^{-})\big)&\\
\\
\mathcal{S}_{1}(Z)(d,u)=\frac{1}{2}\big(Z(d,T^{+})-Z(d,T^{-})\big)&
\end{array}
\labell{re3.3}\end{equation}
where
$T^{+\frac{1}{2}}=T^{+}(\frac{d}{2},u)$,
$T^{-\frac{1}{2}}=T^{-}(\frac{d}{2},u)$ and $Z(d,u)$ is a  formal
power series in $d$ whose coefficients are polynomials. 
We can show that
\begin{eqnarray}
\begin{array}{ll}
 \mathcal{C}_{1}=2\mathcal{S}^{2}+Id 
\\
  \mathcal{S}_{1}=2\mathcal{S}\mathcal{C}
\end{array}
\labell{re4}\end{eqnarray}
and 
\begin{eqnarray}
\begin{array}{llll}
\mathcal{C}_{1}(Q\cdot G)=\mathcal{C}_{1}(Q)\mathcal{C}_{1}(G)+\mathcal{S}_{1}(Q)\mathcal{S}_{1}(G)&\\
\\
\mathcal{S}_{1}(Q\cdot G)=\mathcal{C}_{1}(Q)\mathcal{S}_{1}(G)+\mathcal{S}_{1}(Q)\mathcal{C}_{1}(G))&\\
\\
    \mathcal{C}(Q\cdot G)=\mathcal{C}(Q)\mathcal{C}(G)+\mathcal{S}(Q)\mathcal{S}(G)&\\
\\
    \mathcal{S}(Q\cdot G)=\mathcal{C}(Q)\mathcal{S}(G)+\mathcal{S}(Q)\mathcal{C}(G)&
\end{array}
\labell{re5}\end{eqnarray}
if $Q,G$ are formal power series in $d$ whose coefficients are polynomials of $u$.

\section{Norms for polynomials and basis }
In this section we recall some definitions and results of \cite{SV}. Using a certain suquence of polynomials. we define convenient norms on  spaces of polynomials which satisfies some useful proprieties. We denote by 
\begin{itemize}
\item $\mathcal{P}$ the set of all polynomial whose coefficents are complex
\item $\mathcal{P}_{n}$ the spaces of all polynomials of degree less
  than or equal to n 
\end{itemize}

\begin{propo} \cite{SV}. We  define the sequence of polynomials $\tau_{n}(u)$ by 
$$\tau_0(u)=1,\,\tau_{1}(u)=u,\,\tau_{n+1}(u)=\frac{1}{n}D\tau_{n}(u) \mbox{ for }n\geq1,$$ where the operator $D$ is defined by $$D:=(1-u^{2})\frac{\partial}{\partial{u}}.$$
Then we have
\begin{enumerate}
\item $T^{+}(d,u)=\sum^{\infty}_{n=0}\tau_{n+1}(u)d^{n}$,
%\item $\tau_{1}(u)=u, \tau_{n+1}(u)=\frac{1}{n}D\tau_{n}(u)$
\item  $\tau_n(u)$ has exactly degree $n$ and hence $\tau_0(u),...,\tau_{n}(u)$ form a basis of $\mathcal{P}_{n}$,

\item $\tau_{n}(\tanh(z))=\frac{1}{(n-1)\,!}\big(\frac{d}{dz} \big)^{n-1}\big(\tanh(z)\big)$
\end{enumerate}

\labell{prop1}\end{propo}
\begin{defi} Let $p \in \mathcal{P}_{n}$.  As $
  \tau_0(u),...,\tau_{n}(u)$ form a basis of $\mathcal{P}_{n}$, we can write
  $p\in {\mathcal P}_n$ as 
\begin{eqnarray*}
 p=\sum^{n}_{k=0}a_{k}\tau_{k}(u).
\end{eqnarray*}
Then we define the norms
 \begin{equation}
\|p\|_{n}=\sum^{n}_{i=0}|a_{i}|\left(\frac{\pi}{2}\right)^{n-i} .
\labell{norme}\end{equation}
\end{defi}\
\begin{theo}\cite{SV}. Let n,m be positive integers and $p \in \mathcal{P}_{n}$, 
  \,$q \in \mathcal{P}_{m}.$
  The norms (\ref{norme}) have the following properties:
\begin{enumerate}

\item  $ \|Dp\|_{n+1}\le n\|p\|_{n}$.   
\item If the constant term of $p$ in the basis $\{\tau_0,\tau_1..,\tau_n\}$ is zero,  we have $$\|p\|_{n}\le \|Dp\|_{n+1}.$$
\item There exists a constant $M_{2}$ such that $\|pq\|_{n+m}\le M_{2} \|p\|_{n}\|q\|_{m}$.
\item There is a constant $M_{3}$ such that that for all $ n> 1$,  $|p(u)|\le
  M_{3}\left(\frac{2}{\pi}\right)^{n} \|p\|_{n}$ $(-1\le u\le 1)$.

\item  There is a constant $M_{4}$ such that for all $ n> 1$  with
  $p(1)= p(-1)=0$
\begin{eqnarray*}
\Big\|\frac{p}{\tau_{2}}\Big\|_{n-2}\le M_{4} \|p\|_{n}
\end{eqnarray*}
\end{enumerate}
\labell{th2}\end{theo}
\section{Operators}
In this section we will use some definitions of Schäfke-Volkmer\cite{SV} and adapt their results on operators on polynomial series to our context. 
Let $$\mathcal{Q}:=\left\{\mathnormal{Q}(d,u)=\sum^{\infty}_{n=0}\mathnormal{Q}_{n}(u)d^{n},\textrm{ where}\  \mathnormal{Q}_{n}(u)\in \mathcal{P}_{n},\textrm{ for all}\  n\in \mathbb{N} \right\}.$$
By abuse of notation, let $\|Q\|_n=\|Q_n\|_n$ for a polynomial series $$\mathnormal{Q}(d,u)=\sum^{\infty}_{n=0}\mathnormal{Q}_{n}(u)d^n.$$
\begin{defi}
Let $f$ be  formal power series of $z$ whose coefficients are complex. 
We define a linear operator $f(dD)$ on $\mathcal{Q}$ by
\begin{equation}
f(dD)\mathnormal{Q}(d,u)=\sum^{\infty}_{n=0} \Big( \sum^{n}_{i=0} f_{i}
  D^{i}\mathnormal{Q}_{n-i}(u) \Big)d^{n}
\end{equation}
where $f(z)=\sum^{\infty}_{i=0}f_{i}z^{i}$ and $ \mathnormal{Q}\in\mathcal{Q}$.
\end{defi}
\noindent By the above Definition and (1) of  Proposition  \ref{prop1}  we can show that 
\begin{eqnarray*} 
\mathnormal{Q}(d,T^{+}\big(\theta d, u)\big)=\big(\exp(\theta
dD)\mathnormal{Q}\big)(d,u) \ \textrm{for }\ \mathnormal{Q}\in \mathcal{Q}
 \  \textrm{and all}\  \theta\in \mathbb{C}
\end{eqnarray*} 
Thus with (\ref{re4}) and (1) of  Proposition  \ref{prop1} we obtain
\begin{equation}          
\begin{array}{llll}
    \mathcal{C}(\mathnormal{Q})\
    =\cosh(\frac{d}{2}D)\mathnormal{Q},\ \mathcal{S}(\mathnormal{Q})\ =\sinh(\frac{d}{2}D)\mathnormal{Q}&\\ 
\mathcal{C}_{1}(\mathnormal{Q})=\cosh(dD)\mathnormal{Q},\
\mathcal{S}_{1}(\mathnormal{Q})=\sinh(dD)\mathnormal{Q}
\end{array}
\end{equation}
for polynomial series  $\mathnormal{Q}$ in $ \mathcal{Q}$. 

{\bf Remark.}
According  to the  definition of  norms in    (\ref{norme}), we  have
\begin{eqnarray}
\text{If} \ \ \ Q\in \mathcal{Q},\ \ \  \text{then}   \ \ \ dQ \in \mathcal{Q} \ \ \ \text{and} \ \ \  \| dQ\|_{n}=\dfrac{\pi}{2}\| Q\|_{n-1}\ \  \text{for all}\ \ n\geq 1.
\labell{remp}\end{eqnarray}
\begin{theo}\cite{SV}
Let $f(z)$ be  formal power series having a radius of convergence
greater than $2\pi$ and let k be a positve integer. There is a constant 
$K$ such that: 
If $\mathnormal{Q}$ is a polynomial series having the following property 

\begin{eqnarray*} 
\|\mathnormal{Q}\|_{n}\leq \left\{\begin{array}{ll} 0 & \textrm{for}\  n<k
      \\ M(n-k)\,!(2\pi)^{-n}  & \textrm{for}\  n\geq k
 \end{array} \right.
\end{eqnarray*} 
where $M$ is independent of $n$ and  $\mathnormal{Q} \in  \mathcal{Q}$
then the polynomial series $f(dD)\mathnormal{Q}$ satisfies 
\begin{eqnarray*} 
\|f(dD)\mathnormal{Q}\|_{n}\leq \left\{\begin{array}{ll} 0 & \textrm{for}\  n<k
      \\ MK(n-k)\,!(2\pi)^{-n}  & \textrm{for}\  n\geq k
 \end{array} \right.
\end{eqnarray*} 
\labell{th3}\end{theo}
 Now we define on $\mathcal{Q}$ the following operator
\begin{equation}
\mathcal{J}=\frac{\mathcal{S}}{dD}.
\labell{re6}\end{equation}
where the notation $\frac{\mathcal{S}}{dD}$ means simply $F(dD)$ with  $F(z)=\dfrac{1}{z}\sinh(\frac z2)$.
\begin{lem}
For  each  integer $k$  there exist a positive constant $K$ such that:
If $\mathnormal{Q}$ is a polynomial series 
with odd $\mathnormal{Q}_{n}$ of degree at most $n$, $\|\mathnormal{Q}\|_{n}=0$ for
$n<k$  in case of positive $k$ and 
\begin{eqnarray*} 
\|dD\mathnormal{Q}\|_{n}\leq M(n-k)\,!(2\pi)^{-n}  & \textrm{for}\   n\geq \max(0,k),
\end{eqnarray*} 
 where $M$ is independent of $n$, 
 then the  polynomial series
  $\mathcal{J}^{-1}(\mathnormal{Q})$ satisfies
\begin{eqnarray*} 
\quad \|\mathcal{J}^{-1}(\mathnormal{Q})\|_{n}\leq  MK(2\pi)^{-n}
\left\{\begin{array}{ll}
%(n+3)\,! & \textrm{for}\  k=-2
(n-k+1)\,! & \textrm{for}\  k\leq1
%\\ n\,! & \textrm{for}\  k=1
 \\ (n-1)\,!\log(n)  & \textrm{for}\   k=2 
\\ (n-1)\,! & \textrm{for}\ k \geq 3 
 \end{array} \right.
\end{eqnarray*}
\labell{lem4}\end{lem}
\noindent \textbf{Proof}. We can see easily  that $\mathcal{J}^{-1}=\pi \,\tilde{\mathcal C}^{-1}+g(dD)$, where $\tilde{\mathcal C}=\cosh(\frac14dD)$ and 
$g(z)$ is analytic for  $|z|<4\pi$, and use the proof of \cite{SV}.

We have ${\mathcal S}=dD\,{\mathcal J}={\mathcal J}\, dD$, but using this relation for the inversion of $\mathcal S$ would give an insufficient result. Using of the formula
$$1=\frac2z\sinh(\frac z2)+F(z)z, \mbox{ where }F(z)=z^{-2}(z-2\sinh(\frac z2))$$
is an entire function, we obtain the relation
\begin{equation}
Q=2{\mathcal J}Q+F(dD)\,dDQ\labell{SJ}
\end{equation}
 for polynomial series $Q\in{\mathcal Q}$.
This will be essential in the proof of 
\begin{theo}
For  each  integer $k$  there exist a positive constant $K$ such that:
If $\mathnormal{Q}$ is a polynomial series 
with odd $\mathnormal{Q}_{n}$ of degree at most $n$, $\|\mathnormal{Q}\|_{n}=0$ for
$n<k$ in case of positive $k$, and 
\begin{eqnarray*} 
\|\mathcal{S}(\mathnormal{Q})\|_{n}\leq M(n-k)\,!(2\pi)^{-n}  & \textrm{for}\   
n\geq \max(0,k),
\end{eqnarray*} 
 where $M$ is independent of $n$, 
 then the  polynomial series
  $\mathnormal{Q}$ satisfies
%\begin{eqnarray*} 
%(i)-\quad \|\mathcal{C}^{-1}(\mathnormal{Q})\|_{n}\leq  MK(2\pi)^{-n} \left\{\begin{array}{ll}
% n\,! & \textrm{for}\  k=1
% \\ (n-1)\,!\log(n)  & \textrm{for}\   k=2 
%\\ (n-1)\,! & \textrm{for}\ k \geq 3 
 %\end{array} \right.
%\end{eqnarray*} 
\begin{eqnarray*} 
\quad \|\mathnormal{Q}\|_{n}\leq  MK(2\pi)^{-n}
\left\{\begin{array}{ll}
%(n+3)\,! & \textrm{for}\  k=-2
\\(n-k+1)\,! & \textrm{for}\  k\leq1
%\\ n\,! & \textrm{for}\  k=1
 \\ (n-1)\,!\log(n)  & \textrm{for}\   k=2 
\\ (n-1)\,! & \textrm{for}\ k \geq 3 
 \end{array} \right.
\end{eqnarray*}
\labell{th4a}\end{theo}
\textbf{Proof.} By the preceding theorem, we have the wanted inequalities for
$dDQ={\mathcal J}^{-1}{\mathcal S}Q$ in the place of $Q$. Here we used again
$\|dDZ\|_{n}\leq (n-1)\|Z\|_{n-1}$ for any polynomial series $Z\in{\mathcal Q}$.
Using theorem \ref{th3} implies the same for $F(dD)dDQ$ with the entire function
$F$ of (\ref{SJ})
As $\|Z\|_n\leq\|dDZ\|_{n+1}$ by theorem \ref{th2}, we find the wanted inequalities
(and even something better in the cases $k\geq2$) also for ${\mathcal J}Q$ because
$dD{\mathcal J}={\mathcal S}$. 
Thus formula (\ref{SJ}) yields the result $\square$

In order to obtain an asymptotic approximation for the coefficients of the formal solution, we will need to reverse some operators. This is not possible for the operators $\mathcal{S}$ and $dD$ on the set $\mathcal{Q}$, but we can define a subset $\mathcal{Q}^{*}$ of $\mathcal{Q}$  on which these operators have a  right inverses. 

If we define 
\begin{eqnarray*}
\mathcal{Q}^{*}:=\left\{\mathnormal{Q}(d,u)=\sum^{\infty}_{n=1}\mathnormal{P}_{n}(u)d^{n},\textrm{ where}\  \mathnormal{P}_{n}(u)\in \mathcal{P}^{*}_{n},\textrm{ for all}\  n\geq 1 \right\}.
\end{eqnarray*}
where $\mathcal{P}_{n}^{*}$  is the subspace of $\mathcal{P}_{n}$ defined by
\begin{eqnarray*}
\mathcal{P}_{n}^{*}:=\big\{\sum^{n}_{i=0}\alpha_{i}\tau_{i}\in \mathcal{P}_{n},| \  \alpha_0=0\  \big\}
\end{eqnarray*}
Then, the restrictions of the operators $dD,\mathcal{S}$ to  $\mathcal{Q}^{*},$ denoted here by the same symbols 
\begin{eqnarray*}
dD&:& \mathcal{Q}^{*}\to (1-u^2)d^2\mathcal{Q}\\
\mathcal{S} &:& \mathcal{Q}^{*}\to (1-u^2)d^2\mathcal{Q}
\end{eqnarray*}
are bijective. We denote by $\mathcal{T}$ the  inverse of  the restriction of  $\mathcal{S}$ to  $\mathcal{Q}^{*},$ ,  and we have 
$$\mathcal{T}\mathcal{S}=\text{Id}  \ \ \text{on }  \ \ \mathcal{Q}^{*}$$

\begin{theo}\cite{SV}
We consider a polynomial series 
\begin{eqnarray*} 
\mathnormal{Q}_{\alpha}(d,u)=\sum^{\infty}_{\substack{n=2} }\alpha_{n}(n-1)\,!\Big(
\frac{i}{2\pi}\Big)^{n-1}\tau_{n}(u)d^{n}
\end{eqnarray*}
where $\alpha_{n}=\mathcal{O}(n^{-k})$ as $n
\to \infty$ with some integer $k\geq 2$. Let $$\alpha:=
\frac{4}{\pi}\sum^{\infty}_{\substack{n=1} }\alpha_{n}, $$ then the  coefficients  
$\left\{ \mathcal{T}(\mathnormal{Q}_{\alpha}) \right\}_{n}$ of\,
$\mathcal{T}(\mathnormal{Q}_{\alpha}) $ satisfy
\begin{eqnarray*} 
\Bigg\| \left\{ \mathcal{T}(\mathnormal{Q}_{\alpha})
\right\}_{n}-\alpha (n-1)\,!\Big(\frac{i}{2\pi}\Big)^{n-1}\tau_{n}\Bigg\|_{n}=\mathcal{O}\Big((n-k)\,!(2\pi)^{-n}\Big)
\end{eqnarray*}
as $n\to \infty$ for $n.$
\labell{th5}\end{theo}
\noindent \textbf{Proof}. The proof of this theorem is completely analogous  to that of  \cite{SV}.

\begin{theo}\cite{SV}
Let $k,l, p,q$ be  integer with $p\geq k$ and $q\geq l$. Define
$m$ as the minimum of $k+q$ and $l+p$. then there is a constant $K$
with the following property:\\
\indent If $\mathnormal{P}$ and $\mathnormal{Q}$ are polynomial series 
such that $\|\mathnormal{P}\|_{n}=0 $ for
$n<p$,\ $\|\mathnormal{Q}\|_{n}=0 $ for $n<q$ and 
\begin{eqnarray*}
\begin{array}{ll}
 \|\mathnormal{P}\|_{n}\leq M_{1}(n-k)\,!\, (2\pi)^{-n}  & \textrm{for}\,\,   
 n\geq p \\ 
\|\mathnormal{Q}\|_{n}\leq  M_{2}(n-l)\,!\, (2\pi)^{-n}  & \textrm{for}\,\,   
 n\geq q  
\end{array}
\end{eqnarray*}
then
\begin{eqnarray*}
\|\mathnormal{PQ}\|_{n}\leq KM_{1} M_{2}(n-m)\,!\, (2\pi)^{-n}  & \textrm{for}\,\,   
 n\geq p+q.
\end{eqnarray*}
\labell{th6}\end{theo}

%\begin{theo}\cite{SV}
%Let $\mathnormal{Q}_{1}(d,u)$ be a convergent polynomial series which
%is even with respect to both variables and has constant term 1.\\
% Let $\mathnormal{Q}_{2}(d,u)=d^{2}(1-u^{2})\mathnormal{Q}_{1}(d,u)$
% and
% $\mathnormal{P}(d,u)=\mathcal{S}(\mathnormal{Q}_{2})/\mathcal{C}(\mathnormal{Q}_{2})$. Consider the linear operator defined by 
%
%\begin{equation}
%\mathcal{L}(\mathnormal{Q})=\mathcal{S}(\mathnormal{Q})-\mathnormal{P}(d,u)\cdot\mathcal{C}(\mathnormal{Q}),\quad 
%\mathnormal{Q}\in \mathcal{Q}
%\labell{re7}\end{equation}
%Then there exist a constant $K$ with the following property. If
%$\mathnormal{Q}$ is an odd polynomial series with odd coefficients
%$\mathnormal{Q}_{n}(u)$ satisfying $\mathnormal{Q}_{n}(1)=0$ for
%all\ $n$, $\mathcal{L}(\mathnormal{Q})=0$ for $n<6$ and 
%\begin{eqnarray*}
%\|\mathcal{L}(\mathnormal{Q})\|_{n}\leq M (n-6)\,!(2\pi)^{-n}\  \textrm{for}\ 
%n\geq 6,
%\end{eqnarray*}
%with $M$ independence of $n$, then also 
%\begin{eqnarray*}
%\|dD\mathnormal{Q}\|_{n}\leq K M (n-6)\,!(2\pi)^{-n}\  \textrm{for}\ 
%n\geq 6.
%\end{eqnarray*}
%\labell{th7}
%\end{theo}

\begin{rem} \labell{remM} \rm Observe that the results of this section can also
be applied, if
the constants $M$ are replaced by any increasing sequence $(M_n)_{n\in N}$.
In theorems \ref{th3} and \ref{th6} the
first $n$ terms of the resulting polynomial series only depend of the first
$n$ terms of the given series, so the "$M$" in the result simply has to be replaced by
"$M_n$".
In lemma \ref{lem4} and theorem \ref{th4a}, the
first $n$ terms of the result depend of the first
$n+1$ given terms, so  "$M$" in the result has to be replaced by
"$M_{n+1}$".\end{rem}

\section{Asymptotic approximation of the coefficients of the formal solution} 
 In this    section we will estimate the coefficients of the formal
 solution obtained previously (section 2). 
 The idea is to write equation (\ref{eq1}) essentially in the form

\begin{equation}
V(d,u)\mathcal{S}(Q_1\mathcal{S}(Q_2\, A))(d,u)=g\big(d,u,A(d,u)\big),
\labell{form}\end{equation}
 where $V,Q_1$ and $Q_2$ are known   polynomials
  of $d$ and  $u$ and $g$ is a certain function of  $d, u$ and $A$ involving the operators  $\mathcal{S}, \mathcal{C}$ and $\mathcal{J}$ multiplied by sufficiently high powers of $d$.

Thanks to this equation, we will %be able to use the inverse the operator $\mathcal{J}=\frac{\mathcal S}{dD}$ and  to
 estimate the coefficients of the formal solution 
 using the results of the previous section. We show that the coefficients of this formal solution is Gevrey-1, more precisely $\|A\|_{n}=O\big(n\,!\,(2\pi)^{-n}\big)$.
 \subsection{Rewriting of equation (\ref{eq1})}
 Consider the decomposition
 %Thanks to this equation we will be able to have a recurrence on $n$,
% using the norms  on the polynomials of orders $n$ defined in (\ref{norme}),
 %which makes it possible to reverse the operator $\mathcal{J}$ and to
 %have a formal estimate of the coefficients of the formal solution $A(d, u)$.  For this purpose we consider the decomposition 
%It will turn out to be convenient to consider the decomposition 

\begin{equation}
  A(d,u)=U(d,u)+F(d,u) 
\labell{re8}\end{equation}
where $U$ is the initial part of $A$ calculated before
\begin{equation*}
  U(d,u)=-\frac{1}{4}ud^{2}+\bigg(\frac{91}{864}u^{3}-\frac{47}{576}u \bigg)d^{4}
 +\bigg(-\frac{319}{2880}u^{5}+\frac{185}{1152}u^{3}-\frac{3703}{69120}u \bigg)d^{6}.
\end{equation*}
 We insert this into (\ref{eq2}), with  (\ref{re4}) and (\ref{re5}), and obtain
\begin{equation}
2\cosh(d)\cdot\mathcal{C}_{1}(F)-2u\,\sinh(d)\cdot\mathcal{S}_{1}(F)=W_0\cdot F +f_{1}\Big(d,u,F(d,u)\Big)
\labell{re9}\end{equation}
where 
\begin{eqnarray*}
W_0&=&\Big( \cosh^{2}(d)-u^{2}\sinh^{2}(d)\Big)\bigg[2+(2u^{2}-1)\varepsilon^{2}\cos\Big(U\cdot\sqrt{1-u^{2}}\Big)\\
&-&2u\,\varepsilon^{2}\sin\Big(U\cdot\sqrt{1-u^{2}}\Big)\cdot\sqrt{1-u^{2}}\bigg]\\
f_{1}\Big(d,u,F(d,u)\Big)&=& \varepsilon^{2}\Big( \cosh^{2}(d)-u^{2}\sinh^{2}(d)\Big)\Bigg[\bigg(\dfrac{2u^{2}-1}{\sqrt{1-u^{2}}}\sin\Big(U\cdot\sqrt{1-u^{2}}\Big)\\
&+&2u\cos\Big(U\cdot\sqrt{1-u^{2}}\Big)\bigg)\cos\Big(F(d,u)\sqrt{1-u^{2}}\Big) +\\
&&\bigg(\dfrac{2u^{2}-1}{\sqrt{1-u^{2}}}\cos\Big(U\cdot\sqrt{1-u^{2}}\Big)-2u\sin\Big(U\cdot\sqrt{1-u^{2}}\Big)\bigg)\times\\
&&\bigg(\sin\Big(F(d,u)\sqrt{1-u^{2}}\Big)-F(d,u)\sqrt{1-u^{2}}\bigg)\Bigg]\\
&-&\big(\cosh^{2}(d)-u^{2}\sinh^{2}(d)\big)\sum_{n=1}^{+\infty}I_{2n-1}(u)d^{2n}-2\cosh(d)\mathcal{C}_{1}(U)\\
&-&2u\,\sinh(d)\mathcal{S}_{1}(U)-2\Big( \cosh^{2}(d)-u^{2}\sinh^{2}(d)\Big)\cdot U.
\end{eqnarray*}
Observ that $f_{1}$ has the form  
\begin{equation}
\begin{array}{rcl}f_{1}\big(d,u,F(d,u)\big)&=&\ds y_0(d,u)+y_1(d,u)\sum_{n=1}^{\infty}
\frac1{(2n)!}(1-u^2)^nF(d,u)^{2n}+\\&&\ \ \ \ \ \ \ \ \ \ds y_2(d,u)\sum_{n=1}^{\infty}
\frac1{(2n+1)!}(1-u^2)^nF(d,u)^{2n+1},
\end{array}\labell{ff}\end{equation} 
where  $y_{n}(d,u), n=1,2,3$  are  convergent polynomial series.

Now, we let 
\begin{eqnarray}
\begin{array}{lll}
 F(d,u)&=&Q(d,u)\cdot G(d,u),\\
J(d,u)&=&Q_{1}(d,u)\cdot\mathcal{S}(G),
\end{array}
\labell{re13}\end{eqnarray}
where $ G $ is a formal
power series whose the first term contains $d^{8}$ and 

\begin{eqnarray}
\begin{array}{lll}
Q(d,u)&=&1+\frac{1}{4}(1-u^{2})d^{2}+\bigg(\frac{91}{432}u^{4}-\frac{13}{48}u^{2}+\frac{13}{216}\bigg)d^{4}\\
&+&\bigg(-\frac{319}{960}u^{6}+\frac{1079}{1728}u^{4}-\frac{937}{2880}u^{2}+\frac{287}{8640}\bigg)d^{6},\\
Q_{1}(d,u)&=&(u^{2}-1)d^{2}+\frac{1}{4}(1-u^{4})d^{4}-\dfrac{5}{48}(1-u^{2})\bigg(\frac{4}{9}u^{4}+u^{2}+1\bigg)d^{6}\\
&+&(1-u^{2})\bigg(-\frac{367}{2160}u^{6}+\frac{185}{432}u^{4}-\dfrac{997}{4320}u^{2}\bigg)d^{8}.
\end{array}
\labell{re10}\end{eqnarray}
The choice of $Q_{1}(d,u)$ and  $Q(d,u)$ depends in a precise way of the form of the equation (\ref{form}) and has been determined using Maple.

Using (\ref{re13}), (\ref{re10}) and  (\ref{re5}), we can rewrite equation (\ref{re9}) in the form
%We replace this in  equation (\ref{re9}) and  use  (\ref{re5}).  We obtain 
\begin{eqnarray*}
W_0Q+f_{1}\Big(d,u,F_{d}(u)\Big)&=&\Big[2\cosh(d)\mathcal{C}_{1}(Q)-2u\,\sinh(d)\mathcal{S}_{1}(Q)\Big]\mathcal{C}_{1}(G)\nonumber\\
&+&\Big[2\cosh(d)\mathcal{S}_{1}(Q)-2u\,\sinh(d)\mathcal{C}_{1}(Q)\Big]\mathcal{S}_{1}(G).
\end{eqnarray*}
 Using (\ref{re4}), we obtain 

\begin{eqnarray}
V\cdot \mathcal{S}^{2}(G)+W\cdot \mathcal{S}\mathcal{C}(G)=W_{1}G+f_{2}\Big(d,u,F_{d}(u)\Big)
\labell{re11}\end{eqnarray}
where $f_{2}\Big(d,u,F_{d}(u)\Big)=\frac{1}{4}f_{1}\Big(d,u,F_{d}(u)\Big)$

\begin{eqnarray}
\begin{array}{lll}
V(d,u)&=&\cosh(d)\mathcal{C}_{1}(Q)-u\,\sinh(d)\mathcal{S}_{1}(Q)\\
W(d,u)&=&\cosh(d)\mathcal{S}_{1}(Q)-u\,\sinh(d)\mathcal{C}_{1}(Q)\\
W_{1}(d,u)&=&\frac{1}{4}\Big(-2\cosh(d)\mathcal{C}_{1}(Q)+2u\,\sinh(d)\mathcal{S}_{1}(Q)+W_0Q\Big)
\end{array}
\labell{re12}\end{eqnarray}
The calculation of the first terms of the series $W_{1}$ by Maple shows that the convergent polynomial series  $  W_{1}(d,u)$  begins with a term containing $d^{10}$.

%{\bf RE13}
%\begin{eqnarray}
%\begin{array}{lll}
%J(d,u)&=&Q_{1}(d,u)\cdot\mathcal{S}(G)\\
%V_{1}(d,u)&=&1+(1-u^{2})d^{2}+\bigg(-\dfrac{71}{432}u^{4}-\dfrac{1}{12}u^{2}+\dfrac{107}{432}\bigg)d^{4}\\
%&+&\bigg(\dfrac{1655}{1728}u^{6}-\dfrac{3395}{1728}u^{4}+\dfrac{3289}{2880}u^{2}-\dfrac{389}{2880}\bigg)d^{6}
%\end{array}
%\labell{re13}\end{eqnarray}
%where 

%{\bf RE14}
%\begin{eqnarray}
%Q_{1}(d,u)&=&(u^{2}-1)d^{2}+\frac{1}{4}(1-u^{4})d^{4}+\bigg(\frac{5}{108}u^{6}+\frac{25}{432}u^{4}-\frac{5}{48}\bigg)d^{6}\nonumber\\
%&+&\bigg(\frac{721}{864}u^{8}-\frac{1085}{432}u^{6}+\frac{11035}{4320}u^{4}-\frac{119}{135}\bigg)d^{8}
%\labell{re14}\end{eqnarray}
Using (\ref{re13}) and  (\ref{re5}), we find
\begin{equation}
\mathcal{S}(J)=\mathcal{S}\Big(Q_{1}\mathcal{S}(G)\Big)=\mathcal{C}(Q_{1})\mathcal{S}^{2}(G)+\mathcal{S}(Q_{1})\mathcal{S}\mathcal{C}(G).
\labell{re15}\end{equation}
Using
\begin{eqnarray*}
V_{1}(d,u)&=&1+(1-u^{2})d^{2}+\bigg(-\dfrac{71}{432}u^{4}-\dfrac{1}{12}u^{2}+\dfrac{107}{432}\bigg)d^{4}\\
&+&\bigg(\dfrac{1351}{2160}u^{6}-\dfrac{193}{144}u^{4}+\dfrac{49}{60}u^{2}-\dfrac{11}{108}\bigg)d^{6},
\end{eqnarray*}
we obtain
\begin{eqnarray*}
V_{1}\cdot\mathcal{C}(Q_{1})&=&Q_{1}V+W_{2}\\
V_{1}\cdot\mathcal{S}(Q_{1})&=&Q_{1}W+W_{3}
\end{eqnarray*}
where $W_{2}(d,u) $  and   $ W_{3}(d,u)$ are convergent polynomials series beginning with $d^{10}$. With
(\ref{re11}) and (\ref{re15}), this implies 

\begin{equation}
V_{1}\cdot\mathcal{S}\Big(Q_1\mathcal{S}(G)\Big)=W_{2}\cdot\mathcal{S}^2(G)+W_{3}\mathcal{S}\mathcal{C}(G)+Q_{1}W_{1}G+Q_{1}f_{2}\Big(d,u,Q_1 G(d,u)\Big).
\labell{re16}\end{equation}
This  allows us to prove the  following theorem
\begin{theo}
\begin{eqnarray}
G(d,u)&=&\Big(\frac{\alpha}{d^2}+(\beta+\frac{\alpha}{3})\Big)\dfrac{u\, H_{0}(d,u)}{\tau_2(u)}-\big(\beta d+\frac{\alpha}{d}\big)H_{2}(d,u)\nonumber\\
&+&\delta d H_{1}(d,u)+S(d,u)
\end{eqnarray}
where 
$\alpha, \beta, \delta$ are  constants and the polynomial series $H_0, H_1, H_2, S $  are defined by 
\begin{eqnarray}
H_{0}(d,u):&=&\sum_{\substack{n=10\\ n\  even}}^{\infty}(n-1)\,!\Big(\dfrac{i}{2\pi}\Big)^{n}\tau_{n}(u)\,d^{n}\nonumber\\
H_{1}(d,u):&=&\sum_{\substack{n=9\\ n\  odd}}^{\infty}(n-1)\,!\Big(\dfrac{i}{2\pi}\Big)^{n+1}\tau_{n}(u)\,d^{n}\nonumber\\
H_2(d,u):&=&\sum_{\substack{n=9\\ n\  odd}}^{\infty}n\,!\Big(\dfrac{i}{2\pi}\Big)^{n+1}\tau_{n}(u)\,d^{n}
\labell{rethtt}\end{eqnarray}
and $S(d,u)$ is a polynomial series satisfying $$ \|S\|_{n}=\mathcal{O}\Big( (n-3)! (2\pi)^{n}\Big).$$
\labell{thtt}\end{theo}
To prove this theorem we need to make some overvaluations on the coefficients of the polynomial series $\mathcal{S}\big(Q_1\mathcal{S}(G)\big)$. This will be the subject of the following paragraph.

\noindent{\bf Remark}: Observe that  the series $F,G$ are odd in $u$, even in $d$ and beginning with $d^8$. The series $J$ is even in $u$, odd in $d$ and beginning with $d^{11}$ . In the series $F, G, A, J$,  the degree of the polynomial that is the coefficient of $d^n$ is at most $n-1$
; thus the results of  section 4 can still be applied and  $d^{-1}F, d^{-1}G, d^{-1}A, d^{-1}J \in \mathcal{Q}$%will yield better results if we divide by $d$.

\subsection{Upper bounds for  the coefficients of $\mathcal{S}\big(Q_1\mathcal{S}(G)\big)$}
In this paragraph, we will use  equation (\ref{re16}), together with the definitions of $V_{1}$ and $Q_{1}$, $J$ and $G$, $W_{i}, i=1, 2, 3$, to prove
\begin{lem}
\begin{displaymath}
\big\|\dfrac{1}{d}\mathcal{S}\big(Q_1\mathcal{S}(G)\big)\big\|_{n}=\mathcal{O}\Big( (n-8)\,!(2\pi)^{-n}\Big)
\quad \textrm{as}\  n\to \infty.
\end{displaymath}
\labell{lem8}\end{lem}  

\noindent\textbf{Proof.} We set
\begin{equation}
e_{n}:=\frac{(2\pi)^{n}\big\|\mathcal{S}(J)\big\|_{n}}{ (n-8)\,!}\ \ 
\textrm{for}\ \  n \geq 12 
\labell{re17}\end{equation}
We must show that $e_{n}=\mathcal{O}(n^{-1})$. In the
sequel, we will use the following  convention: if  $a_{n}, n=0, 1, ...$. is any
sequence of positve real numbers, then 
\begin{equation*}
a^{+}_{n}:=\max(a_0,a_{1},...a_{n}) \ \textrm{for all} \ n\geq 0
\end{equation*}
We have 

\begin{equation}
\big\|\mathcal{S}(J)\big\|_{n}\leq e_{n}(n-8)\,! (2\pi)^{-n}\quad\textrm{for}\ \  n \geq 12 .
\labell{re18}\end{equation}
 %The definition of the operator $ \mathcal{J} $,  gives  $ \mathcal{S}=\sinh(dD)=\mathcal{J}(dD)dD$.
Using  Theorem \ref{th4a} and Remark \ref{remM}, we obtain 

\begin{equation}
\big\|J\big\|_{n}\leq K_{1}e^{+}_{n+1}(n-1)\,! (2\pi)^{-n}, \quad\textrm{for}\ \  n \geq 11 .
\labell{re19}\end{equation}
where $K_{1}$ denotes the constant associated with the operator $\mathcal{S}$ in Theorem \ref{th4a}, it is independent  of the present context. In this proof  $K_{i}, i=1,.. 9$  will always denote constants independent of $n$ and the sequence $e_{n}$.
%according to the theorem \ref{th2}
%{\bf RE20}
%\begin{equation}
%\big\|J\big\|_{n}\leq \big\|dDJ\big\|_{n+1}
%\labell{re20}\end{equation}
%this gives 
%{\bf RE21}
%\begin{equation}
%\big\|J\big\|_{n}\leq e^{+}_{n+1}n\,! (2\pi)^{-n-1}
%\labell{re21}\end{equation}

Using (\ref{re13}), we obtain
\begin{equation}
\big\|Q_{1}\mathcal{S}(G)\big\|_{n}\leq K_{1} e^{+}_{n+1}(n-1)\,! (2\pi)^{-n} \quad\textrm{for}\ \  n \geq 11
\labell{re22}\end{equation}
We use (5) of  Theorem \ref{th2} and \ref{th6}. Since $
\dfrac{Q_{1}(d,u)}{\tau_{2}(u)d^{2}} $ is a convergent power series beginning with
$1$, there  is a constant $K_{2}$ %independent of $n$ and the sequence $e_{n}$ 
such that 
\begin{equation}
\big\|\mathcal{S}(G)\big\|_{n}\leq K_{2}\,e^{+}_{n+3}(n+1)\,! (2\pi)^{-n}, \quad\textrm{for}\ \  n \geq 9.
\labell{re23}\end{equation}
Using again  Theorem \ref{th4a} (and remark \ref{remM}) and the fact that $F=G/Q$ where
$Q$ is given in (\ref{re10}), we  obtain 
%{\bf RE24}
%\begin{equation}
%\big\|dDG\big\|_{n}\leq M_{1}e^{+}_{n+3}(n+3)\,! (2\pi)^{-n}
%\labell{re24}\end{equation}
%therefore
\begin{equation}\renewcommand{\arraystretch}{1.3}
\begin{array}l
\big\|G\big\|_{n}\leq K_{3}e^{+}_{n+4}(n+2)\,! (2\pi)^{-n}\quad\textrm{for}\ \  n \geq 8,
\\
\big\|F\big\|_{n}\leq K_{3}e^{+}_{n+4}(n+2)\,! (2\pi)^{-n}\quad\textrm{for}\ \  n \geq 8,
\end{array}
\labell{re25}\end{equation}
where $K_{3}$ is a constant independent of $n$.

This together with theorem \ref{th6} implies that there are constants $K_4,L$ such that
for all $k\geq2$
\begin{equation}
\big\|F^{k}\big\|_{n} \leq K_{4}L_1^kf^{(k)}_{n}(n-5)\,!(2\pi)^{-n}\ \ \ \textrm{for}\ 
n\geq 8k
\labell{re26}\end{equation}
where  
\begin{displaymath}\renewcommand{\arraystretch}{1.4}\begin{array}{lcl}
f^{(2)}_{n}&=&\ds
 \sum^{n-8}_{i=8}e^{+}_{i+4}e^{+}_{n-i+4} \frac{(i+2)\,!\,
  (n-i+2)\,!}{(n-5)\,!},\ \ \textrm{for\ $n\geq 16$},\\
f^{(k+1)}_n&=&\ds
 \sum^{n-8k}_{i=8}e^{+}_{i+4}f^{(k)}_{n-i} \frac{(i+2)\,!\,
  (n-i-5)\,!}{(n-5)\,!},\ \ \textrm{for\ $n\geq 8(k+1)$},
\end{array}
\end{displaymath}
with $f^{(k)}_{n}:=0$ for $n<8k$.

Using Theorems \ref{th3} and \ref{th6} and $W_{i}=O(d^{10}), i=1, 2, 3$, we obtain
\begin{equation}
\Big\|W_{2}\cdot\mathcal{S}^2(G)\big\|_{n}\leq K_{5}e^{+}_{n-6}(n-9)\,! (2\pi)^{-n}\ \ \ \textrm{for}\ n\geq 20
\labell{re26.1}\end{equation}
\begin{equation}
\Big\|W_{3}\mathcal{S}\mathcal{C}(G)\Big\|_{n}\leq K_{6}e^{+}_{n-6}(n-9)\,! (2\pi)^{-n}\ \ \ \textrm{for}\ n\geq 19,
\labell{re26.2}\end{equation}
\begin{equation}
\Big\|Q_{1}W_{1}G\Big\|_{n}\leq K_{7}e^{+}_{n-8}(n-10)\,! (2\pi)^{-n}\ \ \ \textrm{for}\ n\geq 20,
\labell{re26.3}\end{equation}
\begin{equation}
\Big\|Q_{1}f_{2}\Big(d,u,F_{d}(u)\Big)\Big\|_{n}\leq K_{8}
(1+\sum_{k\geq2}\frac{L_2^k}{k!}f^{(k)+}_{n-4})(n-9)\,! (2\pi)^{-n} \ \ \ \textrm{for}\ n\geq 12,
\labell{re26.4}\end{equation}
%with $K_{i}, i=1,.. 4$ independent of $n$ and the sequence $e_{n}$.

Now, let us take  the equation (\ref{re16})
\begin{equation*}
\Big\|V_{1}\cdot\mathcal{S}(J)\Big\|_{n}\leq\Big\|W_{2}\cdot\mathcal{S}^2(G)\big\|_{n}+\Big\|W_{3}\mathcal{S}\mathcal{C}(G)\Big\|_{n}+\Big\|Q_{1}W_{1}G\Big\|_{n}+\Big\|Q_{1}f_{2}\Big(d,u,F_{d}(u)\Big)\Big\|_{n}
\end{equation*}
Using (\ref{re26.1}), (\ref{re26.2}), (\ref{re26.3})  and (\ref{re26.4}), we obtain 
\begin{equation}
\big\|V_{1}\cdot\mathcal{S}(J)\big\|_{n}\leq K_{9}\Big(1+e_{n-6}^{+}+
\sum_{k\geq2}\frac{L_2^k}{k!}f^{(k)+}_{n-4}\Big)(n-9)!\,(2\pi)^{-n}
\end{equation}
Since, $ V_{1} $ is a convergent polynomial series begins with $1$ , we also have 
\begin{equation}
\big\|\mathcal{S}(J)\big\|_{n}\leq K_{10}\Big(1+e_{n-6}^{+}+\sum_{k\geq2}\frac{L_2^k}{k!}f^{(k)+}_{n-4}\Big)(n-9)!\,(2\pi)^{-n} 
\end{equation}

Using (\ref{re17}), we obtain 
\begin{equation}
e^{+}_{n}\leq \dfrac{K}{n}\Big(1+e_{n-6}^{+}+\sum_{k\geq2}\frac{L_2^k}{k!}f^{(k)+}_{n-4}
\Big) \ \ \ \textrm{for}\ n\geq 12
\labell{re27}\end{equation}

{\bf lem5.2}
\begin{lem} Under the condition (\ref{re27}),  we have $e_{n}=\mathcal{O}(n^{-1})$ as
$n\to \infty .$
\labell{lem5.2}\end{lem} 
\noindent\textbf{Proof.} 
Let $K_{1}\geq 10\,!e^{+}_{12}$ an arbitrary number. We assume that 
\begin{equation}
e_{n}\leq \frac{K_{1}(n+p-1)\,!}{(n-2)\,!(p+11)\,!}\ \ \qquad  \textrm{for $12\leq 
  n\leq N-4 $}
\labell{re27.1}\end{equation}
with some $p\geq -1,\, N\geq 16$. This gives for $16\leq  n\leq N $
\begin{eqnarray*}
(n-5)\,!f^{(2)}_{n}\leq K_{1}^{2}\sum_{i=8}^{n-8}\frac{(i+p+3)\,!(n+p-i+3)\,!}{\big((p+11)\,!\big)^{2}}.
\end{eqnarray*} 
The first and last term of the above sum are the largest, so we can easily estimate
\begin{eqnarray*}
\sum^{n-8}_{i=8}(i+p+3)\,!(n+p-i+3)\,!\leq (p+11)\,!(n+p-4)\,!\ \ \  .
\end{eqnarray*} 
We obtain 
\begin{eqnarray*}
f^{(2)}_{n}\leq K^2_{1}\frac{(n+p-4)\,!}{(p+11)\,!(n-5)\,!}  
\qquad  \textrm{for $16\leq  n\leq N $} . 
\end{eqnarray*} 
In a similar way, we can prove by induction that
\begin{eqnarray*}
f^{(k)}_{n}\leq K^k_{1}\frac{(n+p-4)\,!}{(p+11)\,!(n-5)\,!}  
\leq K^k_{1}\frac{(n+p-1)\,!}{(p+11)\,!(n-2)\,!}\qquad  \textrm{for $8k\leq  n\leq N $} . 
\end{eqnarray*}
Using the assumption of the lemma, we obtain 
\begin{eqnarray*}
e_{n}\leq\frac{K}{n}e^{(1+L_2)K_1}\frac{(n+p-1)\,!}{(p+11)\,!(n-2)\,!}  \qquad  \textrm{for $16\leq  n\leq N $}.
\end{eqnarray*} 

Now we choose $N_0\geq 16$ so large that $\frac{K\exp((1+L_2)K_1)}{N_0}\leq
K_{1}$ and then $p$ so large that (\ref{re27.1}) holds for $N=N_0$.
In a first step, our considerations imply by induction over $N$ that
(\ref{re27.1}) holds for {\em all} $N$ and hence 
\begin{eqnarray*}
e_{n}=\mathcal{O}\bigg(\frac{(n+p-1)\,!}{(n-2)\,!}\bigg) \ \textrm{as $  n\to \infty $}
\end{eqnarray*}
for this possibly large value of $p$.

As $K_1$ is arbitrary in (\ref{re27.1}), we have also shown for any $p\geq -1$ that
\begin{eqnarray*}
e_{n}=\mathcal{O}\bigg(\frac{(n+p-1)\,!}{(n-2)\,!}\bigg) \ \textrm{as $  n\to \infty $}
\end{eqnarray*}
implies that
\begin{eqnarray*}
e_{n}=\mathcal{O}\bigg(\frac{(n+p-2)\,!}{(n-2)\,!}\bigg) \ \textrm{as $
  n\to \infty $\ \ .}
\end{eqnarray*}
Consequently the last assertion is proved for $p=-1$ and we have shown
\begin{eqnarray*}
e_{n}=\mathcal{O}\big(n^{-1}\big) \ \textrm{as $
  n\to \infty $ }\ \ .
\end{eqnarray*}
Finally  we have proved that
\begin{eqnarray*}
\big\|\mathcal{S}(J)\big\|_{n} =\mathcal{O}\big((n-9)\,!(2\pi)^{-n}\big) \ \textrm{as $ n\to \infty $ } 
\end{eqnarray*}
and hence that 
\begin{eqnarray*}
\big\|\dfrac{1}{d}\mathcal{S}(J)\big\|_{n} =\mathcal{O}\big((n-8)\,!(2\pi)^{-n}\big) \ \textrm{as $ n\to \infty $ } 
\end{eqnarray*}
 which completes the proof of the lemma. \quad $\square$
\subsection{Proof of theorem \ref{thtt}}

Let $E:=\dfrac{1}{d}\mathcal{S}(J)=\mathcal{S}(d^{-1}J)$. The polynomial series $E$ is odd in 
$d$ and its coefficients are odd in $u$. We partition it
\begin{eqnarray}
E_{n}(u)&=&\alpha_{n}(n-1)\,!\Big(\frac{i}{2\pi}\Big)^{n-1}\tau_{n}(u)+\beta_{n-2}(n-3)\,!\Big(\frac{i}{2\pi}\Big)^{n-3}\tau_{n-2}(u)+\nonumber\\
&&\gamma_{n-4}(n-5)\,!\Big(\frac{i}{2\pi}\Big)^{n-5}\tau_{n-4}(u)+\overline{E}_{n-6}(u)
\labell{re28}\end{eqnarray}
for odd $n\geq 11$, where $\alpha_{n}$ ,  $\beta_{n}$  and  $\gamma_{n}$   are real
numbers and also $\overline{E}_{n}$ have at most degree $n$ for all $n$.
For the whole series $E$ this is equivalent to 
\begin{eqnarray}
\mathcal{S}(d^{-1}J)=E=E_{1}+d^{2}E_{2}+d^{4}E_{3}+d^{6}\overline{E}
\labell{re29}\end{eqnarray}
where 
\begin{eqnarray*}
E_{1} &=&\sum^{+\infty}_{n=11}\alpha_{n}(n-1)\,!\Big(\frac{i}{2\pi}\Big)^{n-1}\tau_{n}(u)d^{n} \\ 
E_{2}&=&\sum^{+\infty}_{n=9}\beta_{n}(n-1)\,!\Big(\frac{i}{2\pi}\Big)^{n-1}\tau_{n}(u)d^{n}\\ 
E_{3}&=&\sum^{+\infty}_{n=7}\gamma_{n}(n-1)\,!\Big(\frac{i}{2\pi}\Big)^{n-1}\tau_{n}(u)d^{n}\\ 
\overline{E} &=&\sum^{+\infty}_{n=5}\overline{E}_{n}(u)d^{n} 
\end{eqnarray*}
Lemma \ref{lem8} implies that
\begin{eqnarray*}
\alpha_{n}=\mathcal{O}(n^{-7}),
\ \beta_{n}=\mathcal{O}(n^{-5}),\ \gamma_{n}=\mathcal{O}(n^{-3}) 
\ \textrm{and} \ \|\overline{E}_{n}\|_{n}=\mathcal{O}\big((n-2)\,!(2\pi)^{-n} \big).
\end{eqnarray*}
Applying $ \mathcal{T}$ to (\ref{re29}) we obtain
\begin{equation}
\frac1dJ=\mathcal{T}(E_{1})+d^{2}\mathcal{T}(E_{2})+d^{4}\mathcal{T}(E_{3})+d^{6}\mathcal{T}(\overline{E})
\labell{re30}\end{equation}
To the first three summands we apply Theorem \ref{th5}. % applies together with the observation
%that $ \mathcal{S}^{-1}={\mathcal J}^{-1} (dD)^{-1}$ and that $(dD)^{-1}$ operates
%essentially like a shift operator on series like $E_j,j=1,2,3$. 
Thus we obtain
\begin{eqnarray*}
\Big\| \{\mathcal{T}(E_{1})\}_{n}-\alpha
(n-1)\,!\Big(\frac{i}{2\pi}\Big)^{n}\tau_{n}(u)\Big\|_{n}&=&\mathcal{O}\big((n-7)\,!(2\pi)^{-n}
\big)\\
\Big\| \{\mathcal{T}(E_{2})\}_{n}-\beta
(n-1)\,!\Big(\frac{i}{2\pi}\Big)^{n}\tau_{n}(u)\Big\|_{n}&=&\mathcal{O}\big((n-5)\,!(2\pi)^{-n}
\big)\\
\Big\| \{\mathcal{T}(E_{3})\}_{n}-\gamma
(n-1)\,!\Big(\frac{i}{2\pi}\Big)^{n}\tau_{n}(u)\Big\|_{n}&=&\mathcal{O}\big((n-3)\,!(2\pi)^{-n}\big)
\end{eqnarray*}
where 
\begin{eqnarray*}
\alpha=\frac{4}{\pi}\sum^{\infty}_{\substack{n=11\\ n\  odd}}\alpha_{n}, \
\beta=\frac{4}{\pi}\sum^{\infty}_{\substack{n=9\\ n\  odd}}\beta_{n},
\gamma=\frac{4}{\pi}\sum^{\infty}_{\substack{n=7\\ n\  odd}}\gamma_{n}.
\end{eqnarray*}
To the last part of (\ref{re30}) we apply  Theorem \ref{th4a} and obtain
\begin{eqnarray*}
\Big\| \{\mathcal{T}(\overline{E})\}_{n}\Big\|_{n}&=&\mathcal{O}\big((n-1)\,!(2\pi)^{-n}\log(n)
\big),
\end{eqnarray*}
thus altogether
\begin{eqnarray*}
\bigg\| \{(\dfrac{1}{d}J)\}_{n}&-&\alpha
(n-1)\,!\Big(\frac{i}{2\pi}\Big)^{n}\tau_{n}-\beta
(n-3)\,!\Big(\frac{i}{2\pi}\Big)^{n-2}\tau_{n-2}-\\
&&\gamma(n-5)\,!\Big(\frac{i}{2\pi}\Big)^{n-4}\tau_{n-4}\bigg\|_{n}
= \mathcal{O}\Big((n-7)\,!(2\pi)^{-n}\log(n)\Big).
\end{eqnarray*}
%If we use the relation $\tau_{n+1}=\frac{1}{n}D\tau_{n}$ and we apply (2) of   Theorem \ref{th2}, we obtain
Using (\ref{remp}) we obtain 
%or equivalently
\begin{eqnarray*}
\bigg\| \{J\}_{n}&-&\alpha
(n-2)\,!\Big(\frac{i}{2\pi}\Big)^{n-1}\tau_{n-1}-\beta
(n-4)\,!\Big(\frac{i}{2\pi}\Big)^{n-3}\tau_{n-3}\\
&-&\gamma(n-6)\,!\Big(\frac{i}{2\pi}\Big)^{n-5}\tau_{n-5}\bigg\|_{n}
= \mathcal{O}\Big((n-8)\,!(2\pi)^{-n}\log(n)\Big)\ \ .
\end{eqnarray*}
\begin{rem}
The asymptotic of  $J_n$ gives a good approximation of $\alpha$; its  suffice to calculate , using a formal calculation software (for example: Pari), the first $40$ terms of $A(d,u)$ by the recurrence of Section 2 and to evaluate the highest coefficients of $ J_n$ to get the approximation $\alpha= 89.0334.$
\labell{rem1}\end{rem}
Next we observe that $J=Q_{1}\mathcal{S}(G)$, where $Q_1$ is given in (\ref{re10})
%$ Q_{1}=-\tau_{2} d^{2}+ \tau_{2} \big(\frac{1}{2}-\dfrac{\tau_{2}}{4}\big)d^{4}+\mathcal{O}(d^{6}) $.
Using part 5.\  of  Theorem \ref{th2}, we obtain 

{\bf EA}
\begin{eqnarray}
\bigg\| \big\{\mathcal{S}(G)\big\}_{n}&+&\alpha
n\,!\Big(\frac{i}{2\pi}\Big)^{n+1}\frac{\tau_{n+1}}{\tau_{2}}+\big(\beta-\alpha p_{2}(u)\big)
(n-2)\,!\Big(\frac{i}{2\pi}\Big)^{n-1}\frac{\tau_{n-1}}{\tau_{2}}\nonumber\\
&&+(\gamma-\beta  p_{2}(u)-\alpha  p_{4}(u)\big)
(n-4)\,!\Big(\frac{i}{2\pi}\Big)^{n-3}\frac{\tau_{n-3}}{\tau_{2}}\bigg\|_{n}\nonumber\\
&= &\mathcal{O}\Big((n-6)\,! (2\pi)^{-n}\log(n)\Big)
\labell{ea}\end{eqnarray}
%\begin{eqnarray*}
%\bigg\| \big\{\mathcal{S}(G)\big\}_{n}&+&\alpha
%n\,!\Big(\frac{i}{2\pi}\Big)^{n+1}\frac{\tau_{n+1}}{\tau_{2}}+\beta
%(n-2)\,!\Big(\frac{i}{2\pi}\Big)^{n-1}\frac{\tau_{n-1}}{\tau_{2}}\\
%&+&\dfrac{\alpha}{4}
%(n-2)\,!\Big(\frac{i}{2\pi}\Big)^{n-1}(2-\tau_{2})\frac{\tau_{n-1}}{\tau_{2}}\bigg\|_{n}= 
%\mathcal{O}\Big((n-3)\,! (2\pi)^{-n}\log(n)\Big)
%\end{eqnarray*}
where 
\begin{eqnarray*}
p_{2}(u)&=&-\dfrac{1}{4}-\dfrac{u^{2}}{4}=-\dfrac{1}{2}+\dfrac{\tau_{2}}{4}\\
p_{4}(u)&=&\dfrac{1}{24}-\dfrac{u^{2}}{48}-\dfrac{7 u^{4}}{432}=\dfrac{1}{216}+\dfrac{55}{1296}\tau_{2}+\dfrac{7}{432}\tau_{4}.
\end{eqnarray*}
%\begin{eqnarray*}
%\bigg\| \big\{\mathcal{S}(G)\big\}_{n}&+&\alpha
%n\,!\Big(\frac{i}{2\pi}\Big)^{n+1}\frac{\tau_{n+1}}{\tau_{2}}+\big(\beta+\dfrac{\alpha}{2}\big)
%(n-2)\,!\Big(\frac{i}{2\pi}\Big)^{n-1}\frac{\tau_{n-1}}{\tau_{2}}\\
%&-&\dfrac{\alpha}{4}
%(n-2)\,!\Big(\frac{i}{2\pi}\Big)^{n-1}\tau_{n-1}\bigg\|_{n}= 
%\mathcal{O}\Big((n-3)\,! (2\pi)^{-n}\log(n)\Big)
%\end{eqnarray*}

{\bf Remark:} Observe that the approximation (\ref{ea}) of the coefficients  $\{\mathcal{S}(G)\}_n$  is polynomial.
 Indeed;  the polynomials $\tau_n(u),  n\geq 2 $ are divisible by  $\tau_2 (u)$. 
\\

In order to find an asymptotic estimation for the coefficients of the formal solution, we need to apply  the inverse of  operator $\mathcal{S}$ . To this purpose, we show the following lemma
\begin{lem}  If $H_0, H_1, H_2$ are  the polynomial series  defined in (\ref{rethtt}). 
Then
\begin{itemize}
\item1 If we define the operator $\mathcal{C}_{\frac{1}{4}}=\cosh\Big(\dfrac{dD}{4}\Big)$, then the polynomial series $\mathcal{C}_{\frac{1}{4}}(H_0), \mathcal{C}_{\frac{1}{4}}(H_1)$ are converging.
\item2.   the polynomial series  $\mathcal{S}(H_{0}), \mathcal{S}(H_1) $ are converging.
\item3. $\mathcal{C}(H_{0})=-H_{0}+\mu_{1}(d,u) $, where $\mu_{1}(d,u) $ is a convergent series.
\item4. $\mathcal{S}(H_{2})=\dfrac{1}{2}H_{0}+\mu_{2}(d,u) $, where $H_{2}=d\dfrac{\partial}{\partial d} H_{1}$ and  $\mu_{2}(d,u) $ is a convergent series .
\item5.  $\mathcal{S}\bigg(\dfrac{H_{0}}{\tau_{2}(u)}\bigg)=-u\sinh(d)\dfrac{H_{0}}{\tau_{2}(u)}+\mu_{3}(d,u) $, where $\mu_{3}(d,u) $ is a convergent series .
\item6.  $\mathcal{S}\bigg(\dfrac{u\,H_{0}}{\tau_{2}(u)}\bigg)=-\frac{1}{2}(u^{2}+1)\sinh(d)\dfrac{H_{0}}{\tau_{2}(u)}+\mu_{4}(d,u) $, where $\mu_{4}(d,u) $ is a convergent series .
\item7.  $\mathcal{S}\bigg(\dfrac{u\,H_{0}}{\sinh(d)\tau_{2}(u)}- H_{2}\bigg)=-\dfrac{H_{0}}{\tau_{2}(u)}+\mu_{5}(d,u) $, where $\mu_{5}(d,u) $ is a convergent series .
\end{itemize}
\end{lem}
\noindent{\bf Proof.} (1)- We have % define the operators  $\mathcal{C}_{\frac{1}{4}}=\cosh\Big(\dfrac{dD}{4}\Big), \mathcal{S}_{\frac{1}{4}}=\sinh\Big(\dfrac{dD}{4}\Big)$, then 
\begin{eqnarray*}
\mathcal{C}_{\frac{1}{4}}(H_{0})=\sum_{\substack{m=0\\ m\  even}}^{\infty}\dfrac{1}{4^{m}m!}d^mD^{m}\bigg(\sum_{\substack{n=10\\ n\  even}}^{\infty}(n-1)\,!\Big(\dfrac{i}{2\pi}\Big)^{n}\tau_{n}(u)\,d^{n}\bigg)
\end{eqnarray*}
using the definition of the operator $D$ in Proposition \ref{prop1}, we obtain
\begin{eqnarray}
\mathcal{C}_{\frac{1}{4}}(H_{0})&=& \sum_{m, n}\dfrac{1}{4^{m}m!}(n+m-1)\,!\Big(\dfrac{i}{2\pi}\Big)^{n}\tau_{n+m}(u)\,d^{n+m}\nonumber\\
&=& \sum_{\substack{k=10\\ k\  even}}^{\infty}\gamma_{k}(k-1)!\Big(\dfrac{i}{2\pi}\Big)^{k}  \tau_k d^k
\labell{pl}\end{eqnarray}
where
$$\gamma_{k}:= \sum_{\substack{m=0\\ m\  even}}^{k-10}\dfrac{1}{4^{m}m!}\big(\frac{2\pi}{i}\big)^{m}.$$
Hence
\begin{eqnarray*}
\gamma_{k}= \sum_{\substack{m=0\\ m\  even}}^{\infty}\dfrac{1}{4^{m}m!}\big(\frac{2\pi}{i}\big)^{m}-\sum_{\substack{m=k-8\\ m\  even}}^{\infty}\dfrac{1}{4^{m}m!}\big(\frac{2\pi}{i}\big)^{m}.
\end{eqnarray*}
Using $$\sum_{\substack{m=0\\ m\  even}}^{\infty}\dfrac{1}{4^{m}m!}\big(\frac{2\pi}{i}\big)^{m}=\sum_{l=0}^{\infty}\dfrac{(-1)^l}{(2l)!}\big(\frac{\pi}{2}\big)^{2l}=\cos(\frac{\pi}{2})=0,$$
we find
\begin{eqnarray*}
\gamma_{k}= -\sum_{\substack{m=k-8\\ m\  even}}^{\infty}\dfrac{1}{m!}\big(\frac{\pi}{2i}\big)^{m}
\end{eqnarray*}
which implies 
\begin{eqnarray*}
|\gamma_{k}|&\leq & \dfrac{1}{(k-8)!}\Big(\frac{\pi}{2}\Big)^{k-8}\sum_{\substack{m=0\\ m\  even}}^{\infty}\dfrac{1}{m!}\big(\frac{\pi}{2}\big)^{m}\\
&\leq & \dfrac{\cosh(\pi/2)}{(k-8)!}\Big(\frac{\pi}{2}\Big)^{k-8}
\end{eqnarray*}
This with (\ref{pl}) imply that $\mathcal{C}_{\frac{1}{4}}(H_{0})=\mu(d,u)$ is convergent. For $\mathcal{C}_{\frac{1}{4}}(H_1)$, we can use the same method.

 (2)- As  $ \mathcal{S}=2\mathcal{S}_{\frac{1}{4}}\mathcal{C}_{\frac{1}{4}}$, where $\mathcal{S}_{\frac{1}{4}}=\sinh\Big(\dfrac{dD}{4}\Big)$,  we obtain using (1),
$$\mathcal{S}(H_0)=2\mathcal{S}_{\frac{1}{4}}\Big(\mathcal{C}_{\frac{1}{4}}(H_{0})\Big)=2\mathcal{S}_{\frac{1}{4}}(\mu)$$
%$$\mathcal{S}(H_1)=2\mathcal{S}_{\frac{1}{4}}\Big(\mathcal{C}_{\frac{1}{4}}(H_{1})\Big)=2\mathcal{S}_{\frac{1}{4}}(\mu)$$
This implies that $\mathcal{S}(H_0)$ is convergent. 
 For $\mathcal{S}(H_1)$, we can use the same method.

(3)-  We have $\mathcal{C}(H_{0})=(2\mathcal{C}^{2}_{\frac{1}{4}}-Id)H_{0}=-H_0+2\mathcal{C}^{2}_{\frac{1}{4}}(H_0)$. This with (1) imply  $$\mathcal{C}(H_{0})=-H_{0}+\mu_{1}(d,u)$$ where $\mu_{1}(d,u)$ is  a convergent series.

(4)- We  differentiate the equation $\mathcal{S}(H_{1})=C_{1}$ with respect  to $d$. As $$z\dfrac{d}{dz}\Big(\sinh(z/2)\Big)=\dfrac{z}{2}\cosh(z/2), $$  we obtain
\begin{eqnarray*}
\mathcal{S}(H_{2})+\dfrac{1}{2}\mathcal{C}(dDH_{1})=d\dfrac{\partial C_{1}}{\partial d}
\end{eqnarray*}
Using (3) of this lemma, we obtain 
\begin{eqnarray*}
\mathcal{S}(H_{2})=\dfrac{1}{2}dDH_{1}+\mu_{2}(d,u)=\dfrac{1}{2}H_{0}+\mu_{2}(d,u),
\end{eqnarray*}
 where $\mu_{2}(d,u)$  is a convergent series.

(5)- Using (\ref{re5}) and (2), (3) of this lemma we obtain 
\begin{eqnarray*}
\mathcal{S}(H_{0})=\mathcal{S}\bigg(\tau_{2}(u)\dfrac{H_{0}}{\tau_{2}(u)}\bigg)&=&\mathcal{S}\big(\tau_{2}(u)\big)\mathcal{C}\bigg(\dfrac{H_{0}}{\tau_2}\bigg)+\mathcal{C}(\tau_{2}(u))\mathcal{S}\bigg(\dfrac{H_{0}}{\tau_2}\bigg)\\
&=& \mu_{3}(d,u),\\
\mathcal{C}(H_{0})=\mathcal{C}\bigg(\tau_{2}(u)\dfrac{H_{0}}{\tau_{2}(u)}\bigg)&=&\mathcal{C}\big(\tau_{2}(u)\big)\mathcal{C}\bigg(\dfrac{H_{0}}{\tau_2}\bigg)+\mathcal{S}(\tau_{2}(u))\mathcal{S}\bigg(\dfrac{H_{0}}{\tau_2}\bigg)\\
&=& -H_0.
\end{eqnarray*}
This implis
\begin{eqnarray*}
\mathcal{S}\bigg(\dfrac{H_{0}}{\tau_2}\bigg)&=&\dfrac{\mathcal{S}\big(\tau_{2}(u)\big)}{\mathcal{C}\big(\tau_{2}(u)\big)^2-\mathcal{S}\big(\tau_{2}(u)\big)^2}H_0+\mu_{4}(d,u)\\
&=&-u\sinh(d)\dfrac{H_{0}}{\tau_{2}(u)}+\mu_{4}(d,u),
\end{eqnarray*}
 where $\mu_{3}(d,u)$, $\mu_{4}(d,u)$  are  converging series.

(6)- The proof of (6) is similar to that of (5).

(7)- Using (4) and (6), we obtain 
\begin{eqnarray*}
\mathcal{S}\bigg(\dfrac{u\,H_{0}}{\sinh(d)\tau_{2}(u)}- H_{2}\bigg)&=&\Big(-\frac{1}{2}(u^2+1)-\frac{1}{2}\tau_{2}(u)\Big)\dfrac{H_{0}}{\tau_{2}(u)}+\mu_{5}(d,u)\\
  &=&-\dfrac{H_{0}}{\tau_{2}(u)}+\mu_{5}(d,u),
\end{eqnarray*}
where $\mu_{5}(d,u)$  is a convergent series. This completes the proof of Lemma.
\\

Using the definition of $H_{0}$ in (\ref{rethtt}) and (\ref{ea}), we can rewrite $\dfrac{1}{d}\mathcal{S}(G)$ in the form
\begin{eqnarray}
\frac{1}{d}\mathcal{S}(G)&=&-\alpha \dfrac{H_{0}}{d^{2}\tau_{2}(u)}-(\beta +\frac{1}{2}\alpha)\dfrac{H_{0}}{\tau_{2}(u)}+\dfrac{1}{4}\alpha H_{0}+\overline{X},
\labell{re31}\end{eqnarray}
where 
$$\big\|\overline{X}\big\|_{n}= \mathcal{O}\Big((n-3)\,!(2\pi)^{-n}\Big).$$
Using (2) and (5) of the previous Lemma and applying also the inverse of the operator $\mathcal{S}$ in (\ref{re31}), we obtain 

{\bf RE34}
\begin{eqnarray}
\frac{1}{d}G=\dfrac{1}{d^2}\Big(\alpha+(\beta +\frac{1}{2}\alpha)d^2\Big)\Big[\dfrac{1}{\sinh(d)}\dfrac{u\, H_{0}}{\tau_{2}(u)}-H_{2}\Big]+\dfrac{1}{2}\alpha H_{2}+\delta H_{1}+\overline{X}_{1}
\labell{re34}\end{eqnarray}
where, $\delta$ is a constant and 
\begin{eqnarray}
H_{2}(d,u)&=&d\dfrac{\partial H_{1}}{\partial d}(d,u) =\sum_{\substack{n=9\\ n\  odd}}^{\infty}n\,!\Big(\dfrac{i}{2\pi}\Big)^{n+1}\tau_{n}(u)\,d^{n}\\
\big\|\overline{X}_{1}\big\|_{n}&=& \mathcal{O}\Big((n-2)\,!(2\pi)^{-n}\Big)\\
\delta &=&\dfrac{4}{\pi} \sum_{n=1}^{\infty}\delta_n,
\end{eqnarray}
with $\delta_n=O(n^{-2}).$
  In the expression of $\frac{1}{d}G$, the term $\delta H_1$ comes from the fact that the series  $X$ can be written
$$X(d,u)=D_1(d,u)+d^2 D_2(d,u)$$ where 
\begin{eqnarray*}
D_1(d,u)&=& \sum_{n=8}^{\infty}\delta_{n} (n-1)!\bigg(\dfrac{i}{2\pi}\bigg)^n \tau_n(u) d^n\\
\big\|D_{2}\big\|_{n}&=& \mathcal{O}\Big((n-1)\,!(2\pi)^{-n}\Big)
\end{eqnarray*}
if we  apply  theorem  \ref{th5} on  the series $D_1(d,u)$  and Theorem \ref{th4a} on $D_1(d,u)$, the term $\delta H_1$  appears in the  expression of $\frac{1}{d}G.$
%We calculated $\mathcal{S}^{-1}(X_{1}), \mathcal{S}^{-1}(X_{2})$, and  $\mathcal{S}^{-1}(X_{3})$ (cf Appendix to the next version of this manuscript). Using also Corollary \ref{th4a}, we obtain 

Since $\dfrac{1}{\text{sinh}(d)}=d^{-1}-\dfrac{d}{6}+\mathcal{O}(d^{3})$, we obtain
%Therefore
\begin{eqnarray}
G(d,u)&=&\Big(\frac{\alpha}{d^2}+(\beta+\frac{\alpha}{3})\Big)\dfrac{u\, H_{0}(d,u)}{\tau_2(u)}-\big(\beta d+\frac{\alpha}{d}\big)H_{2}(d,u)\nonumber\\
&+&\delta d H_{1}(d,u)+S(d,u)
\labell{re36}\end{eqnarray}
where  $ \|S\|_{n}=\mathcal{O}\Big( (n-3)! (2\pi)^{n}\Big)$ $\square$
%$ S(d,u)=\sum_{n=8}^{\infty}S_{n}(u)d^n$ and 
%\begin{eqnarray}
%\alpha_{2}(\varepsilon)&=&\Big(\frac{\alpha}{d^2}+(\beta+\frac{\alpha}{3})\Big]\nonumber\\
%\beta_{2}(\varepsilon)&=& \frac{\alpha}{4}+\big(\frac{3\alpha}{8}+\frac{\beta}{4}\big)d^2\labell{re36a}\\
%\gamma_{2}(\varepsilon)&=&\frac{9\alpha}{4}-\big(\frac{\alpha}{16\pi^{2}}+\frac{5\beta}{4}\big)d^{2}\ \ \ .\nonumber
%\end{eqnarray}

Observe that in $d^{-2}\dfrac{u\, H_{0}}{\tau_2(u)},\, d^{-1}H_{2}$ the degree of the coefficients of $d^n$ exceeds $n$.
This is due to the fact that the expressions $u\frac{\tau_{n+2}}{\tau_2}-\tau_{n+1}$
etc., which are of degree $n-1$, were split.

It is not necessary (but would not be difficult) to write down asymptotic 
approximations for the coefficients of $F$, because equations (\ref{re13})  and 
(\ref{re9}) can be used. This completes the proof of the theorem \ref{thtt}  $\square$

\section{Functions and quasi-solutions}
So far, we have shown that equation (\ref{eq1}) has a formal solution and
we have found an asymptotic approximation of the coefficients of the formal
solution. We will use this to construct a quasi-solution, i.e.\ a function that satisfies
equation (\ref{eq1}) except for some exponentially small error. To 
that purpose, we define the functions
\begin{eqnarray}
H_{n}(u):&=&(n-1)\,!\Big(\frac{i}{2\pi}\Big)^{n}\tau_{n}(u)
%H^{2}_{n}(u):&=&(n+1)\,!\Big(\frac{i}{2\pi}\Big)^{n+2}\tau_{n+1}(u)\nonumber\\
%H^{1}_{n}(u):&=&n\,!\Big(\frac{i}{2\pi}\Big)^{n+1}\tau_{n}(u)
\labell{re37}\end{eqnarray} 
and %for $i=1..3$
\begin{eqnarray}
h_{0}(t,u):&=&\sum^{\infty}_{\substack{n=10\\ n\ even}}H_{n}(u)\frac{t^{n-1}}{(n-1)\,!}\\
h_{1}(t,u):&=&\dfrac{i}{2\pi}\sum^{\infty}_{\substack{n=9\\ n\ odd}}H_{n}(u)\frac{t^{n-1}}{(n-1)\,!}\\
h_{2}(t,u):&=&\dfrac{i}{2\pi}\sum^{\infty}_{\substack{n=9\\ n\ odd}}n\,H_{n}(u)\frac{t^{n-1}}{(n-1)\,!}
\labell{re38}\ \ .\end{eqnarray} 
This means that
\begin{eqnarray}
h_{0}(t,u)&=&\frac{i}{2\pi}\sum^{\infty}_{\substack{n=1\\ n\ odd}}\Big(\frac{i}{2\pi}\Big)^{n}\tau_{n+1}(u)\, t^{n}-\Big(H_2\, t+H_4\frac{t^3}{3!} +H_6\, \frac{t^5}{5!}+H_8\, \frac{t^7}{7!}\Big)\nonumber\\
h_{1}(t,u)&=&\big(\frac{i}{2\pi}\big)^{2}\sum^{\infty}_{\substack{n=0\\ n\ even}}\Big(\frac{i}{2\pi}\Big)^{n}\tau_{n+1}(u)t^{n}-\frac{i}{2\pi}\Big(H_1+H_3\frac{t^2}{2!} +H_5\, \frac{t^4}{4!}7H_7\, \frac{t^6}{6!}\Big)
%h_{2}(t,u)&=&\big(\frac{i}{2\pi}\big)^{2}\sum^{\infty}_{\substack{n=2\\ n\ even}}(n+1)\Big(\frac{i}{2\pi}\Big)^{n}\tau_{n+1}(u)t^{n}\\
%&=& h_{1}(t,u)+\big(\frac{i}{2\pi}\big)^{2}\sum^{\infty}_{\substack{n=2\\ n\ even}}n\,\Big(\frac{i}{2\pi}\Big)^{n}\tau_{n+1}(u)t^{n}
\labell{re39}\end{eqnarray}
Using part 4.\ of the  proposition (\ref{prop1}), we obtain
\begin{eqnarray*}
\frac{i}{2\pi}\sum^{\infty}_{\substack{n=1\\ n\ odd}}\Big(\frac{i}{2\pi}\Big)^{n}\tau_{n+1}(u)\, t^{n}=\frac{i}{2\pi}\sum^{\infty}_{\substack{n=1\\ n\ odd}}\frac{1}{n\,!}\Big(\frac{it}{2\pi}\Big)^{n}\frac{d^{n}}{d^{n}\xi}\Big(\tanh(\xi)\Big)
\end{eqnarray*}
%or equivalently
%\begin{eqnarray*}
%h_{1}(t,u)=\frac{-i\,u}{8\pi^{3}(1-u^{2})}\sum^{\infty}_{\substack{n=1\\ n\ odd}}\frac{1}{n\,!}\Big(\frac{it}{2\pi}\Big)^{n}\frac{d^{n}}{d^{n}\xi}\Big(g_{1}(\xi)\Big),
%\end{eqnarray*}
%where \ $g_{1}(\xi)=-2(\tanh(\xi)-\tanh(\xi)^{3})$ and $\xi=\xi(u)=artanh(u)$\\
%and hence that 
%\begin{eqnarray*}
%h(t,u)=\frac{-1}{\pi^{2}(1-u^{2})}\sum^{\infty}_{\substack{n=0\\ n\ even}}\frac{1}{n\,!}\Big(\frac{it}{\pi}\Big)^{n}\frac{d^{n}}{d^{n}\xi}\Big(g(\xi)\Big)-\frac{2u}{\pi^{2}}
%\end{eqnarray*}
This is obviously the difference of two Taylor expansion and thus we can write
\begin{eqnarray}
h_{0}(t,u)&=&\frac{i}{4\pi}\Big[\tanh\big(\xi+\frac{it}{2\pi}\big)-\tanh\big(\xi-\frac{it}{2\pi}\big)\Big]\nonumber\\
&-&\Big(H_2\, t+H_4\frac{t^3}{3!} +H_6\, \frac{t^5}{5!}+H_8\, \frac{t^7}{7!}\Big).
\labell{re40}\end{eqnarray}
Similarly,
\begin{eqnarray*}
h_{1}(t,u)&=&-\frac{1}{4\pi^2}\sum^{\infty}_{\substack{n=0\\ n\ even}}\frac{1}{n\,!}\Big(\frac{it}{2\pi}\Big)^{n}\frac{d^{n}}{d^{n}\xi}\Big(\tanh(\xi)\Big)\nonumber\\
&-&\frac{i}{2\pi}\Big(H_1+H_3\frac{t^2}{2!} +H_5\, \frac{t^4}{4!}7H_7\, \frac{t^6}{6!}\Big)
%&+& 2\sum^{\infty}_{\substack{n=2\\ n\ even}}\frac{1}{(n-1)\,!}\Big(\frac{i}{2\pi}\Big)^{n+2}\frac{d^{n-1}}{d^{n-1}\xi}\Big(\frac{d}{d\xi}\big(\tanh(\xi)\big)\Big)t^{n-1}
\end{eqnarray*}
or equivalently
%\begin{eqnarray*}
%h_{2}(t,u)&=&\frac{t}{16\pi^{4}}\sum^{\infty}_{\substack{n=0\\ n\ even}}\frac{1}{n\,!}\Big(\frac{it}{2\pi}\Big)^{n}\frac{d^{n}}{d^{n}\xi}\Big(g_{1}(\xi)\Big)\\
%&-& \frac{i}{4\pi^{3}}\sum^{\infty}_{\substack{n=1\\ n\ odd}}\frac{1}{n\,!}\Big(\frac{it}{2\pi}\Big)^{n}\frac{d^{n}}{d^{n}\xi}\Big(g_{2}(\xi)\Big),
%\end{eqnarray*}
%where \ $g_{2}(\xi)=1-\tanh(\xi)^{2}$. 
%We rewrite 
\begin{eqnarray}
h_{1}(t,u)&=&\frac{-1}{8\pi^{2}}\Big[\tanh\big(\xi+\frac{it}{2\pi}\big)+\tanh\big(\xi-\frac{it}{2\pi}\big)\Big]\nonumber\\
&-&\frac{i}{2\pi}\Big(H_1+H_3\frac{t^2}{2!} +H_5\, \frac{t^4}{4!}7H_7\, \frac{t^6}{6!}\Big).
\labell{re41}\end{eqnarray}
%Using the same method for  $h_{2}(t,u)$,  we obtain
 %\begin{eqnarray}
%h_{2}(t,u)&=&h_{1}(t,u)-\frac{i\,t}{8\pi^{3}}\sum^{\infty}_{\substack{n=1\\ n\ odd}}\frac{1}{n\,!}\Big(\frac{it}{2\pi}\Big)^{n}\frac{d^{n}}{d^{n}\xi}\Big(g(\xi)\Big)\nonumber\\
%&=&h_{1}(t,u)- \frac{i\, t}{16\pi^{3}}\Big[g_{1}\big(\xi+\frac{it}{2\pi}\big)-g_{1}\big(\xi-\frac{it}{2\pi}\big)\Big].
%\labell{re42}\end{eqnarray}
%where $g_{1}(z)=1-\tanh(z)^2$.

The functional equations for the trigonometric and hyperbolic functions imply that
\begin{eqnarray}
h_{0}(t,u)&=&-\dfrac{1}{2\pi}\frac{(1-u^2)\,\sin\big(\frac{t}{2\pi}\big)\cos\big(\frac{t}{2\pi}\big)}{\cos\big(\frac{t}{2\pi}\big)^{2}+u^{2}\sin\big(\frac{t}{2\pi}\big)^{2}}\nonumber\\
 &+&\dfrac{\tau_2}{4\pi^2}t-\dfrac{\tau_4}{16\pi^4}t^3+\dfrac{\tau_6}{64\pi^6}t^5-\dfrac{\tau_8}{256\pi^8}t^7\nonumber\\
h_{1}(t,u)&=&-\dfrac{1}{4\pi^2}\frac{(u-u^3)\sin\big(\frac{t}{2\pi}\big)^{2}}{\cos\big(\frac{t}{2\pi}\big)^{2}+u^{2}\sin\big(\frac{t}{2\pi}\big)^{2}}\nonumber\\
&+&\dfrac{\tau_1}{4\pi^2}-\dfrac{\tau_3}{16\pi^4}t^2+\dfrac{\tau_5}{64\pi^6}t^4-\dfrac{\tau_7}{256\pi^8}t^6
%&-&\frac{(u-u^3)\,\sin\big(\frac{t}{2\pi}\big)\cos\big(\frac{t}{2\pi}\big)}{2\pi^{3}\Big(\cos\big(\frac{t}{2\pi}\big)^{2}+u^{2}\sin\big(\frac{t}{2\pi}\big)^{2}\Big)^{2}}\nonumber\\
%h_{2}(t,u)&=&h_{1}(t,u)-\dfrac{t}{4\pi^{3}}\frac{(u-u^3)\,\sin\big(\frac{t}{2\pi}\big)\cos\big(\frac{t}{2\pi}\big)}{\Big(\cos\big(\frac{t}{2\pi}\big)^{2}+u^{2}\sin\big(\frac{t}{2\pi}\big)^{2}\Big)^{2}}
\labell{re43}\end{eqnarray}
For fixed real $u$ the functions $h_{k}(.,u), k=0, 1, 2 $ are analytic in $|t|<\rho$, where $\rho=\pi^2$. In the subsequent definition, 
we consider real values  of $u$,   $0< u\leq 1$, here
$h_{k}, k=0, 1, 2$ are also analytic with respect to $t$ on the positive real axis.

We define the functions $\mathcal{H}_{k}(d,u), k=0..3$,  by
\begin{eqnarray}
\mathcal{H}_{0}(d,u):&=&\int^{+\infty}_0e^{-\frac{t}{d}}h_{0}(t,u) d\,t \ \mbox{for }0<u\leq 1 \nonumber\\
\mathcal{H}_{1}(d,u):&=& \int^{+\infty}_0e^{-\frac{t}{d}}h_{1}(t,u) d\,t \ \mbox{for }0<u\leq 1  \nonumber\\
\mathcal{H}_{2}(d,u):&=& \int^{+\infty}_0 e^{-\frac{t}{d}}h_{2}(t,u) d\,t \ \mbox{for }0<u\leq 1. 
%&=&d\dfrac{\partial \mathcal{H}_{1} }{\partial d}(d,u)%+\int^{+\infty}_0 e^{-\frac{t}{d}}h_{3}(t,u) d\,t \ \mbox{for }0<u\leq1
\labell{re44}\end{eqnarray}
%where $$h_{3}(t,u):=- \frac{i\, t}{16\pi^{3}}\Big[g_{1}\big(\xi+\frac{it}{2\pi}\big)-g_{1}\big(\xi-\frac{it}{2\pi}\big)\Big].$$
We have  $\mathcal{H}_{2}(d,u)=d\dfrac{\partial \mathcal{H}_{1} }{\partial d}(d,u)$. Indeed
\begin{eqnarray*}
d\dfrac{\partial \mathcal{H}_{1} }{\partial d}(d,u)&=& \int^{+\infty}_0\Big(\frac{1}{d} e^{-\frac{t}{d}}\Big)t\cdot h_{1}(t,u)  d\,t \\
&=&- \int^{+\infty}_0\dfrac{\partial}{\partial t}\Big(e^{-\frac{t}{d}}\Big)\cdot t\, \cdot h_{1}(t,u)  d\,t\\
&=&\int^{+\infty}_0 e^{-\frac{t}{d}}\big(h_{1}(t,u)+t\cdot \dfrac{\partial}{\partial t}h_{1}(t,u)\big) d\,t\\
&=& \int^{+\infty}_0 e^{-\frac{t}{d}}h_{2}(t,u) d\,t
\end{eqnarray*}

The functions $\mathcal{H}_k(d,.)$ are real analytic; they can be continued analytically
to the interval $-1<u\leq1$ in the following way. Choose some positive number $M$ and
let $\Gamma_1$ the path consisting of the segment from $0$ to $Mi$ and
of the ray $t\mapsto t+Mi, t\geq0$. Let $\Gamma_2$ the symmetric path that could 
also be obtained using $-M$ instead of $M$. Recalling (\ref{re40}), we can also define
\begin{eqnarray}
\mathcal{H}_0(d,u):&=&\frac{i}{4\pi}
\left[\int_{\Gamma_2}e^{-\frac{t}{d}}
\tanh\big(\xi+\frac{it}{2\pi}\big)\,dt-
\int_{\Gamma_1}e^{-\frac{t}{d}}\tanh\big(\xi-\frac{it}{2\pi}\big)\,dt
\right]\nonumber\\
&+&\mu_{0}(d,u), 
\labell{re44a}\end{eqnarray}
where $$\mu_{0}(d,u):=\dfrac{1}{4\pi^2} \tau_2(u) d^2-\dfrac{3}{8\pi^4}\tau_4(u) d^4+\dfrac{15}{8\pi^6}\tau_6(u) d^6-\dfrac{315}{16\pi^8}\tau_8(u) d^8,$$ 
for $-\tanh(\frac2\pi M)<u\leq1$, where $\xi=\mbox{artanh}(u)$, 
because the singularities of $\tanh$ are $i(\frac\pi2+n\pi)$, 
$n$ integer. As $M$ is arbitrary, this defines the analytic continuation
of $\mathcal{H}_0(d,.)$ for $-1<u\leq1$.
Similarly, the real analytic continuations of $\mathcal{H}_k, k=1, 2$ are defined.

In the sequel, we use the operator $\mathcal{C}, \mathcal{S}$ also for functions.
\begin{lem}  Consider the functions $\mathcal{H}_k(d,u), k=0..2$, defined  in (\ref{re44}). 
Then, for $-1<u\leq1$
\begin{itemize}
\item1 For $k=0,1$, 
\begin{eqnarray}
\mathcal{H}_{k}(d,T^{\pm\frac{1}{2}})=-\mathcal{H}_{k}(d,u)+\mu^{\pm}_k(d,u),
\labell{rel}\end{eqnarray}
 where $T^{+\frac{1}{2}}, T^{-\frac{1}{2}}$ are defined in (\ref{re3.3}) and the functions
$\mu^{\pm}_k(d,u), k=0, 1,$ are analytic,  beginning with $d^{10}$, resp $d^9$
\item2.  For $k=0,1$,    $\mathcal{S}(\mathcal{H}_{k})=\mu_k(d,u) $ , where the functions
$\mu_k(d,u), k=0,1,$ are analytic,  beginning with $d^{11}$,  resp $d^{10}$ .
\item3.  For $k=0,1$,    $\mathcal{C}(\mathcal{H}_{k})=-\mathcal{H}_{k}(d,u)+\lambda_k(d,u) $ , where the functions
$\lambda_k(d,u), k=0,1,$ are analytic,  beginning with $d^{10}$, resp  $d^{9}$.
\item4. $\mathcal{S}(\mathcal{H}_{2})=\dfrac{1}{2}\mathcal{H}_{0}(d,u)+\mu_{2}(d,u) $, where   $\mu_{2}(d,u) $ is a analytic function, beginnings with $d^{10}$ .
%\item5.  $\mathcal{S}\bigg(\dfrac{H_{0}}{\tau_{2}(u)}\bigg)=-u\sinh(d)\dfrac{H_{0}}{\tau_{2}(u)}+\mu_{3}(d,u) $, where $\mu_{3}(d,u) $ is a convergent series .
\item5.  $\mathcal{S}\bigg(\dfrac{u\,\mathcal{H}_{0}}{\tau_{2}(u)}\bigg)=-\frac{1}{2}(u^{2}+1)\sinh(d)\dfrac{\mathcal{H}_{0}}{\tau_{2}(u)}+\mu_{4}(d,u) $, where the function $\mu_{4}(d,u) $ is analytic, beginnings with $d^{11}$  .
%\item7.  $\mathcal{S}\bigg(\dfrac{u\,H_{0}}{\sinh(d)\tau_{2}(u)}- H_{2}\bigg)=-\dfrac{H_{0}}{\tau_{2}(u)}+\mu_{5}(d,u) $, where $\mu_{5}(d,u) $ is a convergent series .
\end{itemize}
\labell{lem6.1}\end{lem}
\noindent{\bf Proof.} (1)- For $k=0$ we replace $u$ by $T^{+\frac{1}{2}}$ in (\ref{re40}). Using  (\ref{re44a}) and $\xi(T^{+\frac{1}{2}})=\xi(u)+\frac{1}{2}d$ we obtain for $0<u\leq1$
\begin{eqnarray}
\mathcal{H}_{0}(d,T^{+\frac{1}{2}})=\int^{+\infty}_0e^{-\frac{t}{d}}h_{0}(t,T^{+\frac{1}{2}}) d\,t=\dfrac{i}{4\pi}\mathcal{I}^{+}+ \mu(d,T^{+\frac{1}{2}})
\labell{re45}\end{eqnarray}
%\begin{eqnarray}
%\mathcal{H}(d,T^{+})=\frac{-i\ T^{+}}{16\pi^{3}\tau_{2}(T^{+})}\mathcal{I}^{+}
%\end{eqnarray}
where 
%\begin{eqnarray*}
%\mathcal{I}^{+}&=&\int^{+\infty}_0e^{-\frac{t}{d}}\Big(g_{1}\big(\xi(T^{+})+\frac{it}{2\pi}\big)+g_{1}\big(\xi(T^{+})-\frac{it}{2\pi}\big)\Big)d\,t  
 %\\
%\mathcal{I}^{-}&=&\int^{+\infty}_0e^{-\frac{t}{d}}\Big(g\big(\xi(T^{-})+\frac{it}{\pi}\big)+g\big(\xi(T^{-})-\frac{it}{\pi}\big)\Big)d\,t  
%\end{eqnarray*}
%Because of \, $\xi(T^{\pm})=\xi\pm d$ where $\xi$ and $u$ are  coupled 
%by $\xi=artanh(u)$ we obtain
%which  imply that 
 \begin{eqnarray*}
\mathcal{I}^{+}=\int^{+\infty}_0e^{-\frac{t}{d}}\tanh\big(\xi+\frac{d}{2}+\frac{it}{2\pi}\big)d\,t -\int^{+\infty}_0e^{-\frac{t}{d}}\tanh\big(\xi-\frac{d}{2}+\frac{it}{2\pi}\big)d\,t
\end{eqnarray*}
If we substitute $t+\pi i\,d $ in the first
part,\ $t-\pi i\,d $ in the second part we obtain 
 \begin{eqnarray*}
\mathcal{I}^{+}=-\int^{+\infty-\pi id}_{-\pi id}e^{-\frac{t}{d}}\tanh\big(\xi+\frac{it}{2\pi}\big)d\,t +\int^{+\infty+\pi id}_{\pi id}e^{-\frac{t}{d}}\tanh\big(\xi-\frac{it}{\pi}\big)d\,t\ \ .
\end{eqnarray*}
Now, we apply Cauchy's theorem 
\begin{eqnarray*}
\mathcal{I}^{+}&=&-\int^{+\infty}_0e^{-\frac{t}{d}}\Big(\tanh\big(\xi+\frac{it}{2\pi}\big)-\tanh\big(\xi-\frac{it}{2\pi}\big)\Big)d\,t\\
&+&\int^{-\pi id}_0e^{-\frac{t}{d}}\tanh\big(\xi+\frac{it}{2\pi}\big)d\,t
 -\int^{\pi id}_0e^{-\frac{t}{d}}\tanh\big(\xi-\frac{it}{2\pi}\big)d\,t\ \ .
\end{eqnarray*}
Substituting $t=-i\,ds $ in the second
part,\ $t= i\,ds $ in the third part, we obtain 
\begin{eqnarray*}
\mathcal{I}^{+}&=&\int^{+\infty}_0e^{-\frac{t}{d}}\Big(\tanh\big(\xi+\frac{it}{2\pi}\big)-\tanh\big(\xi-\frac{it}{2\pi}\big)\Big)d\,t\\
&-&2i\,d\int^{\pi}_0\cos(s) \tanh\big(\xi+d \frac{s}{2\pi}\big)ds.
\end{eqnarray*}
With (\ref{re45}) this implies for $0<u\leq1$ 
\begin{eqnarray}
\mathcal{H}_{0}(d,T^{+})= -\mathcal{H}_{0}(d,u)+\mu_{0}^{+}(d,u)
\end{eqnarray}
where
\begin{eqnarray*}
\mu_{0}^{+}(d,u)=\dfrac{d}{2\pi}\int^{\pi}_0\cos(s) \tanh\big(\xi+d\frac{s}{2\pi}\big)ds+ \mu(d,T^{+\frac{1}{2}}).
\labell{re46}\end{eqnarray*}
By real analytic continuation, this formula is valid for $-1<u\leq1$.  We use the same method for $\mathcal{H}_{0}(d,T^{-\frac{1}{2}}), \mathcal{H}_{1}(d,T^{\pm\frac{1}{2}})$ and obtain for $-1<u\leq1$
%\begin{eqnarray}
%\mathcal{I}^{-}&=&-\int^{+\infty}_0e^{-\frac{t}{d}}\Big(g\big(\xi+\frac{it}{\pi}\big)+g\big(\xi-\frac{it}{\pi}\big)\Big)d\,t-4d(1-u^{2})\mathcal{D}^{-}(d,u)
%\end{eqnarray} 
\begin{eqnarray}
\mathcal{H}_{0}(d,T^{-\frac{1}{2}})&=& -\mathcal{H}_{0}(d,u)+\mu^{-}_{0}(d,u)\nonumber\\
\mathcal{H}_{1}(d,T^{+\frac{1}{2}})&=& -\mathcal{H}_{1}(d,u)+\mu_{1}^{+}(d,u)\nonumber\\
\mathcal{H}_{1}(d,T^{-\frac{1}{2}})&=& -\mathcal{H}_{1}(d,u)+\mu_{1}^{-}(d,u)
%\mathcal{H}_{2}(d,T^{+})&=& \mathcal{H}_{2}(d,u)+ d\,\mathcal{H}_{3}(d,u)+d\mathcal{D}_{1}^{+}(d,u)+d^2\mathcal{D}_{2}^{+}(d,u)\nonumber\\
%\mathcal{H}_{2}(d,T^{-})&=& \mathcal{H}_{2}(d,u)-d\, \mathcal{H}_{3}(d,u)+d\mathcal{D}_{1}^{-}(d,u)+d^2\mathcal{D}_{2}^{-}(d,u),
\labell{re47}\end{eqnarray}
%\int^{\pi}_0\sin(s)\frac{\big(u-\tanh(\frac{ds}{\pi})\big)\big(1-\tanh(\frac{ds}{\pi})^{2}\big)}{\big(1-u\tanh(\frac{ds}{\pi})\big)^{3}}d\,s\\
where 
\begin{eqnarray*}
\mu^{-}_{0}(d,u)&=&\dfrac{d}{2\pi}\int^{\pi}_0\cos(s) \tanh\big(\xi-\frac{d\,s}{2\pi}\big)d\,s+ \mu_{0}(d,T^{-\frac{1}{2}})  \\       
\mu^{+}_{1}(d,u)&=&-\dfrac{d}{4\pi^2}\int^{\pi}_0\sin(s) \tanh\big(\xi+\frac{d\,s}{2\pi}\big)d\,s + \mu_{1}(d,T^{+\frac{1}{2}})\\
\mu^{-}_{1}(d,u)&=&-\dfrac{d}{4\pi^2}\int^{\pi}_0\sin(s) \tanh\big(\xi-\frac{d\,s}{2\pi}\big)d\,s + \mu_{1}(d,T^{-\frac{1}{2}})\\
\mu_{1}(d,u)&=&\frac{1}{4\pi^2}\tau_1(u) d -\frac{1}{8\pi^4}\tau_3(u) d^3+\frac{3}{8\pi^6}\tau_5(u) d^5-\frac{45}{16\pi^8}\tau_7(u) d^7
\end{eqnarray*} 
%From the equations (6.4) and (6.6) we write
%this gives via (2.8)

(2)- Using the definition of the operator $\mathcal{S}$ in (\ref{re3.3}) and  (1) of this Lemma, the result is immediate.

(3)-  The proof of (3) is similar to that of (2).

(4)- For $k=1$, we differentiate (\ref{rel}) with respect to $d$. Because $$\dfrac{\partial T^{\pm\frac{1}{2}}}{\partial d}=\pm\frac{1}{2}\Big(1-(T^{\pm\frac{1}{2}})^2\Big),$$ then
\begin{eqnarray*}
\mathcal{H}_2(d,u)\pm\frac{d}{2}\Big(1-(T^{\pm\frac{1}{2}})^2\Big)\dfrac{\partial \mathcal{H}_{1}}{\partial u}(d,T^{\pm\frac{1}{2}})= -\mathcal{H}_{2}(d,u)+d\mu'^{\pm}_{1}(d,u)
\end{eqnarray*}
implies 
\begin{eqnarray*}
\mathcal{S}\big(\mathcal{H}_2\big)&=&-\frac{d}{2}\mathcal{C}\Big((1-u^2)\dfrac{\partial \mathcal{H}_{1}}{\partial u}\Big)+d\big(\mu'^{+}_{1}(d,u)-\mu'^{-}_{1}(d,u)\big)\\
&=& \frac{d}{2}(1-u^2)\dfrac{\partial \mathcal{H}_{1}}{\partial u}+\mu_{2}(d,u)\\
&=& \frac{1}{2} \mathcal{H}_{0}(d,u)+\mu_{2}(d,u)
\end{eqnarray*}
where $\mu_{2}(d,u)$ is analytic function beginings with $d^{10}$.

(5)- Using (\ref{re5}) and (2), (3) of previous lemma , we obtain 
\begin{eqnarray*}
\mathcal{S}\Big(\dfrac{u}{\tau_2}\mathcal{H}_0\Big)&=&\mathcal{S}\Big(\dfrac{u}{\tau_2}\Big)\mathcal{C}\big(\mathcal{H}_0\big)+\mathcal{C}\Big(\dfrac{u}{\tau_2}\Big)\mathcal{S}\big(\mathcal{H}_0\big)\\
&=& -\frac{1}{2}(u^{2}+1)\sinh(d)\dfrac{\mathcal{H}_{0}}{\tau_{2}(u)}+\mu_{4}(d,u) ,
\end{eqnarray*}
 where the function $\mu_{4}(d,u) $ is analytic, beginnings with $d^{11}$. This completes proof of the Lemma.

 In the sequel we consider $u_0\in]-1,0]$.
\begin{propo}
We have
\begin{enumerate}
\item Uniformly for $u_0\leq u\leq1$, 
\begin{eqnarray}
\mathcal{H}_{0}(d,u) &\sim & \sum^{\infty}_{\substack{n=2\\ n\ even}}(n-1)\,!\Big(\frac{i}{2\pi}\Big)^{n}\tau_{n}(u)d^{n}\ 
\textrm{as} \ d\searrow 0 \nonumber\\
\mathcal{H}_{1}(d,u)&\sim&\sum^{\infty}_{\substack{n=3\\ n\ odd}}(n-1)\,!\Big(\frac{i}{2\pi}\Big)^{n+1}\,\tau_{n}(u)d^{n}\ 
\textrm{as} \ d\searrow 0\nonumber\\
\mathcal{H}_{2}(d,u)&\sim&\sum^{\infty}_{\substack{n=3\\ n\ odd}}n\,!\Big(\frac{i}{2\pi}\Big)^{n+1}\,\tau_{n}(u)d^{n}\ 
\textrm{as} \ d\searrow 0 
\labell{re48}\end{eqnarray} 
%\item For $i=1..3$, \ $0\leq\Big|\int^{t}_0\int^{\sigma}_0\int^{\mu}_0h_{i}(x,u)dx\,d\mu\,d\sigma\Big|\leq Kt^{3}$ for all $t\geq 0$ and\ $0< u\leq 1$ 
\item  For $i=1..3$, $\big|\frac{\partial\mathcal{H}_{i}}{\partial u}(d,u)\big|\leq Kd$\quad  for \
  $u_0< u\leq 1\ $ $(d>0)$
\end{enumerate}
\labell{propo2}\end{propo}
\noindent\textbf{Proof.}
The proof of this proposition is similar to that of  \cite{S}.
%\begin{enumerate}
%\item to prove (2) we use Watson's lemma and (6.2)
%\item follows immediately from
%\end{enumerate}
%\begin{lem}
%We have 
%\begin{eqnarray*}
%\mathcal{H}^{2}(d,u)=\frac{1}{d^{2}}\int_0^{\infty}e^{-\frac{t}{d}}\,f(t,u)dt, 
%\quad (d>0,\ 0<u<1)
%\end{eqnarray*}
%where $f(.,u)$ is analytic in $|t|<\frac{\pi^{2}}{2}$ for $0\leq u\leq 
%1$, $f$ is continuous on $[0,\infty[\times ]0,1]$ and 
%\begin{eqnarray*}
%|f(s,u)|&\leq& Kt^{3} \quad (0\leq t\leq\frac{2\pi^{2}}{3},\ 0<u\leq 1)\\
%\int_{\frac{\pi^{2}}{2}}^{t}e^{-\frac{s}{d}}\,f(s,u)ds &\leq& Kt(t-\frac{\pi^{2}}{2})^{2}, 
%\quad (t\geq\frac{\pi^{2}}{2},\ 0<u\leq 1)\\
%%\end{eqnarray*}
%\end{lem}

With the aim of applying the results of \cite{SV}, we consider $S_{n-2}(u)=R_n(u)$,
where $S_n(u)$ is the remainder term in (\ref{re36}). Then  $R_{n}(u), n$ 
is a  sequence of polynomials of degree at most $n$ and 
\begin{eqnarray*}
\|R_{n}\|_{n}=O\Big((n-5)\,!(2\pi)^{-1}\log(n)\Big)
\end{eqnarray*}
\begin{lem}\cite{SV}
If we define 
\begin{eqnarray*}
r(t,u):&=&\sum_{n=10}^{\infty}R_{n}(u)\frac{t^{n-1}}{(n-1)\,!} \  \quad
(t\in \mathbb{C}, |t|\leq\pi^{2}, u_0\leq u\leq 1 )\\
r(t,u):&=&r\big(\pi^{2},u\big)+\big(t-\pi^{2}\big)\frac{\partial
  r}{\partial t}\big(\pi^{2},u\big) \quad
( t>\pi^{2}, u_0\leq u\leq 1 )\\
\mathcal{R}(d,u):&=&\int_0^{\infty}e^{-\frac{t}{d}}\,r(t,u)dt 
\end{eqnarray*}
then
\begin{enumerate}
\item $r$ is continuously differentiable function on the set $B$ of
all $(t,u)$ such that $u$ satisfies $u_0\leq u\leq 1$ and $t$ is a
complex number and satisfies  $|t|\leq\pi^{2}$ or $
t>\pi^{2}$. The restriction of $r$ to  $u_0\leq u\leq
1,|t|\leq\pi^{2}$ is twice continuously differentiable. for
fixed $u_0\leq u\leq 1$ the function $r(t,u)$ is analytic in
$|t|<\pi^{2}$
\item $\mathcal{R}(d,u)$ is continuous, partially differentiable with
  respect to $u$, has  continuous partial derivative and 
\begin{eqnarray}
\mathcal{R}(d,u)\sim\sum^{\infty}_{n=10}R_{n}(u)d^{n}\ 
\textrm{as} \ d\searrow 0
\end{eqnarray}
\item  $\big|\mathcal{R}(d,u)\big|\leq Kd^{3},\
  \big|\frac{\partial\mathcal{R}}{\partial u}(d,u)\big|\leq Kd^{3}$\quad  for \
  $u_0\leq u\leq 1\ $ $(d>0)$
\end{enumerate}
\labell{lem1}\end{lem}

The importance of our definition of $\mathcal{R}$ lies in a certain
compatibility with insertion of the functions $T^{+}, T^{-}$ for
$u$. First let 
\begin{eqnarray*}
\sum^{\infty}_{n=10}R^{+}_{n}(u)d^{n}&=&\sum^{\infty}_{n=11}R_{n}(T^{+})d^{n}\\
\sum^{\infty}_{n=10}R^{-}_{n}(u)d^{n}&=&\sum^{\infty}_{n=11}R_{n}(T^{-})d^{n}
\end{eqnarray*}
We obtain a new sequences $R_{n}^{+}(u), R_{n}^{-}(u)$ of polynomials
of degree at most $n$. This follows from  the relation 
\begin{equation}
p\big(T^+(d,u)\big)=\sum_{k=0}^{\infty}\frac{1}{k !}D^k p(u) d^k\ \ .
\end{equation}
Theorem \ref{th3}  implies
\begin{eqnarray*}
\big\|R_{n}^{+}(u)\big\|_{n} &=&O\Big((n-5)\,!(2\pi)^{-n}\log(n)\Big)\\
\big\|R_{n}^{-}(u)\big\|_{n}  &=&O\Big((n-5)\,!(2\pi)^{-n}\log(n)\Big)
\end{eqnarray*}
Therefore we can use the previous lemma  for $R_{n}^{+}(u),\  R_{n}^{-}(u)$
and obtain  functions $\mathcal{R}^{+}(d,u),\ \mathcal{R}^{-}(d,u)$.
\begin{theo}
There is a positive constant $K$ independent of $ d,u$ such that 
\begin{eqnarray*}
\big|\mathcal{R}^{+}(d,u)-\mathcal{R}(d,T^{+})\big| &\leq & Kd^{3}e^{-\frac{\pi^{2}}{d}},\quad \text{for}\,   (d>0, u_0<u\leq 1),\\
\big|\mathcal{R}^{-}(d,u)-\mathcal{R}(d,T^{-})\big|&\leq& Kd^{3}e^{-\frac{\pi^{2}}{d}}\quad \text{for}\,      (d>0, u_0<u\leq 1).
\end{eqnarray*}
\labell{th9}\end{theo}
\noindent {\bf Proof.}
\rm The proof is exactly the one of  \cite{SV} .
\begin{defi}

Let $\mathcal{D}(d,u)$ be a function defined for $0<d<d_0$ and
$u_0<u<1$. We say that $\mathcal{D}(d,u)$ has property \textmd{G} if 
\begin{eqnarray*}
\mathcal{D}(d,u)=\int_0^{\infty}e^{-\frac{t}{d}}\,q(t,u)dt \quad
(0<d<d_0,\ u_0<u<1 )
\end{eqnarray*}
is the Laplace transform of some function $q(t,u)$ that has the
following properties  
\begin{enumerate}
\item  $q(t,u)$ is defined if $u_0<u<1$  and either $t$ is complex and
  $|t|<\pi^{2}$ or $t$ is real and $t\geq 0$,
\item $q(t,u)$ is analytic in $|t|<\pi^{2}$ for $ u_0<u<1$,
\item $q(t,u)$ restricted to $0\leq t<\pi^{2}$ or
  $t\geq\pi^{2}$ is continuous and the
  $\lim_{t\to\pi^{2} } q(t,u)$ exists for every  $ u_0<u<1$,
\item there is a positve constant $K$  such that
\begin{eqnarray*}
 |q(d,u)|\leq K e^{Kt} \quad \text{for} \quad t\geq 0, 
(0<d<d_0,\ u_0<u<1 )
\end{eqnarray*}
\end{enumerate}
\end{defi}
%The set of all functions with the property \textmd{G} will be denote by $\mathcal{G}$.
\begin{lem} For $u_0<u\leq 1$, we have

\begin{enumerate}
\item If  $\mathcal{H}_{i}(d,u), i=0,1,2$ are the functions of (\ref{re44}) then 
\begin{equation*}
d^{2}\mathcal{H}_{k}(d,u)=(1-u^2)\tilde{\mathcal{H}}_{k}(d,u)+\mathcal{O}\Big((1-u^2)e^{-\frac{\pi^2}{d}}\Big), \ \ k=0,1
\end{equation*} 
and 
\begin{equation*}
d^{3}\mathcal{H}_{2}(d,u)=(1-u^2)\tilde{\mathcal{H}}_{2}(d,u)+\mathcal{O}\Big((1-u^2)e^{-\frac{\pi^2}{d}}\Big), 
\end{equation*} 
 where $\tilde{\mathcal{H}}_{i}(d,u),i=1,2,3$ have  property \textmd{G}. %Moreover, for every function $\mathcal{D}(d,u) $ with property \textmd{G} the product  $d^{2}\mathcal{D}\cdot\mathcal{H}$ also has  property \textmd{G}
\item  Let $k$ be a psitive integer. If  $\mathcal{D}_{1}, \mathcal{D}_{2}$ have  property \textmd{G} and  their first terms in the Taylor development at $d = 0$, begin with $d^k$  then 
$$\mathcal{D}_{1}(d,u)\mathcal{D}_{2}(d,u)=d^k\mathcal{D}(d,u)+\mathcal{O}\Big(d^ke^{-\frac{\pi^2}{d}}\Big),$$ where $\mathcal{D}(d,u)$ has  property \textmd{G}

\item Any function  $\mathcal{D}(d,u)$ analytic in a neighborhood of $d=0$ has  property \textmd{G} if  $\mathcal{D}(0,u)=0$ for all $u$,
\item If $\mathcal{R}(d,u)$ is defined by lemma \ref{lem1} then 
  $\dfrac{1}{d^{2}}\mathcal{R}(d,u)$  has  property \textmd{G}
\item If  $\mathcal{D}_{1}, \mathcal{D}_{2} $ have  property
  \textmd{G} then so do $\mathcal{D}_{1}+\mathcal{D}_{2},
 \  \mathcal{D}_{1}-\mathcal{D}_{2}$ and\, $\mathcal{D}_{1}\cdot\mathcal{D}_{2} $
\item If  $\mathcal{D}(d,u)$ has  property \textmd{G} then 
\begin{eqnarray}
\big|\mathcal{D}(d,u)\big|\leq Kd  \quad \big(0<d<\frac{1}{K})
\end{eqnarray}
with some constant $K>0$ independent of $u$.
\end{enumerate}
\labell{lem11}\end{lem}

\noindent {\bf Proof}
\begin{enumerate}
\item   For $i=0$,

(i)-If  $u>0$,
 we have 
\begin{eqnarray}
d^2 \mathcal{H}_{0}=(1-u^2) \int_{0}^{\infty} e^{-\frac{t}{d}}g_2(t,u) dt
\labell{fp}\end{eqnarray}
where $$g_2(t,u)=\dfrac{1}{(1-u^2)}\int_{0}^{t}\int_{0}^{\tau}h_0(s,u) ds\, d\tau $$
$g_2(t,u)$ has a logarithmic singularity at $t_k(s)=(2k+1)\pi^2\pm d\dfrac{2\pi\, s }{\varepsilon}i$ for $(k\geq 0, s>0)$. it is analytic in $|t|<\pi^2$ and 
$\lim_{t\to\pi^2}g_2(t,u)$ exists. 

If  we put 
$$\tilde{\mathcal{H}}_{0}(d,u)=\int_{0}^{\infty} e^{-\frac{t}{d}}\tilde{g}_2(t,u) d\,t $$
where $$\tilde{g}_2(t,u)=
\begin{cases}
g_2(t,u),          &\text{if $t\leq \pi^2$}\\
g_2(\pi^2,u),          &\text{if $t\geq \pi^2$}
\end{cases}
$$
then $\tilde{\mathcal{H}}_{0}(d,u)$ has property \textmd{G} and 
$$d^2\mathcal{H}_{0}(d,u)=(1-u^2)\tilde{\mathcal{H}}_{0}(d,u)+\mathcal{O}\Big((1-u^2)e^{-\frac{\pi^2}{d}}\Big). $$
%\noindent Remark: By real analytic continuation the formula (\ref{fp}) is valid for $-1<u_0 <u<1$.

(ii)- For $-u_0 <u<1$, where $0<u_0 <1$,, we have 
\begin{eqnarray*}
\mathcal{H}_{0}(d,u)&=& \int_{0}^{\infty e^{i\varphi}} e^{-\frac{t}{d}}h_0(t,u) d\,t+2\pi i\sum_{k\geq 0}Res\Big(e^{-\frac{t}{d}}h_0(t,u),t_k(s)\Big)\\
&=&\int_{0}^{\infty e^{i\varphi}} e^{-\frac{t}{d}}h_0(t,u) d\,t+\mathcal{O}\Big((1-u^2)e^{-\frac{\pi^2}{d}}\Big)
\end{eqnarray*}
where $\dfrac{\pi}{2}<\varphi<\dfrac{\pi}{4}$. For $0<u<u_0$, this formula coincides with the formula  $$\int_{0}^{\infty} e^{-\frac{t}{d}}h_0(t,u) d\,t $$ and extends it by  real analytic continuation for $-u_0<u<0.$

This implis 
\begin{eqnarray*}
d^2\mathcal{H}_{0}(d,u)&=&(1-u^2)\int_{0}^{\infty e^{i\varphi}} e^{-\frac{t}{d}}g_2(t,u) d\,t+\mathcal{O}\Big((1-u^2)d^2e^{-\frac{\pi^2}{d}}\Big)
\end{eqnarray*}
we obtain 
\begin{eqnarray*}
d^2\mathcal{H}_{0}(d,u)&-&(1-u^2)\tilde{\mathcal{H}}_{0}(d,u)=(1-u^2)\int_{\Gamma} e^{-\frac{t}{d}}g_2(t,u) d\,\Gamma\\
&+&\mathcal{O}\Big((1-u^2)e^{-\frac{\pi^2}{d}}\Big),
\end{eqnarray*}
where $\Gamma$ is the path following the real line from infinity to $(\pi^2,0)$, then along the vertical line from $(\pi^2,0)$ to $\Big(\pi^2,\pi^2\tan(\varphi)\Big)$ and finally along the line $y=\tan(\varphi) x$ from  $\Big(\pi^2,\pi^2\tan(\varphi)\Big)$ to infinity.

Since $g_2(.,u)$ is bounded on $\Gamma$, then 
$$\Big|d^2\mathcal{H}_{0}(d,u)-(1-u^2)\tilde{\mathcal{H}}_{0}(d,u)\Big|\leq K (1-u^2)e^{-\frac{\pi^2}{d}},$$
where $K$ is positive constant. Finally 
$$d^2\mathcal{H}_{0}(d,u)=(1-u^2)\tilde{\mathcal{H}}_{0}(d,u)+\mathcal{O}\Big((1-u^2)e^{-\frac{\pi^2}{d}}\Big) \ \ \ \text{for} \ \ \  (0<u_0<u\leq1)$$

For $i=1$ we can use the same method.

For $i=2$, we use the same method with 
\begin{eqnarray}
d^3 \mathcal{H}_{2}=(1-u^2) \int_{0}^{\infty} e^{-\frac{t}{d}}g_2(t,u) d\,t
\end{eqnarray}
where $$g_2(t,u)=\dfrac{1}{(1-u^2)}\int_{0}^{t}\int_{0}^{\sigma}\int_{0}^{\tau}h_2(s,u) ds\, d\sigma\,d\tau. $$

\item We assume that  $\mathcal{D}_{1}(d,u), \mathcal{D}_{2}(d,u)$ have  property \textmd{G} and their first terms in the Taylor development at $d=0$ begin with $d^k$. Then 
\begin{eqnarray*}
\mathcal{D}_{1}= \int_{0}^{\infty} e^{-\frac{t}{d}}f(t,u) d\,t\\
\mathcal{D}_{2}= \int_{0}^{\infty} e^{-\frac{t}{d}}g(t,u) d\,t
\end{eqnarray*}
where $f(t,u), g(t,u)$ are analytic in $|t|<\pi^2$ and $f(t,u)=O(t^{k-1}), g(t,u)=O(t^{k-1})$. 
\begin{eqnarray}
\mathcal{D}_{1}(d,u)\mathcal{D}_{2}(d,u)= \int_{0}^{\infty} e^{-\frac{t}{d}}(f*g)(t,u) d\,t
\end{eqnarray}
Since 
\begin{eqnarray*}
 h(t,u)&= &(f*g)(t,u) =\int_{0}^{s}f(t,u)g(t-s,u)ds \\
&=&\int_{0}^{t}f(t-s,u)g(s,u)ds\\
&=& \int_{t/2}^{t}f(s,u)g(t-s,u)ds +\int_{t/2}^{t}f(t-s,u)g(s,u)ds.
\end{eqnarray*}
 For $t<\pi^2$, the function $h(t,u)$ is  $k$ times differentiable  with respect to $t$  and 
\begin{eqnarray*}
 h'(t,u)&= & f(t,u)g(0,u)+f(0,u)g(t,u)-f(\frac{t}{2},u)g(\frac{t}{2},u)\\
&+& \int_{t/2}^{t}f(s,u)g'(t-s,u)ds +\int_{t/2}^{t}f'(t-s,u)g(s,u)ds\\
&=&\int_{t/2}^{t}f(s)g'(t-s)ds +\int_{t/2}^{t}f'(t-s,u)g(s,u)ds\\
&-&f(\frac{t}{2},u)g(\frac{t}{2},u)\\
 h^{(k)}(t,u)&= &\int_{t/2}^{t}f(s)g^{(k-1)}(t-s)ds +\int_{t/2}^{t}f^{(k-1)}(t-s,u)g(s,u)ds\\
&- &\sum_{n=0}^{^k-1}f^{(n)}(\frac{t}{2},u)g^{(k-1-n)}(\frac{t}{2},u)
\end{eqnarray*}
Observe that $ h^{(k)}(t,u)$ is contninuous on $[0,2\pi^2[$, it is analytic for $|t|<\pi^2$. If we put 
 $$\tilde{h}(t,u)=
\begin{cases}
h(t,u),          &\text{if $t< \pi^2$}\\
(t-\pi^2)^k.         &\text{if $t\geq \pi^2$}
\end{cases}
$$
then 
$$ \int_{0}^{\infty} e^{-\frac{t}{d}}\tilde{h}^{(k)}(t,u) d\,t,$$  has property \textmd{G} and 
\begin{eqnarray*}
 (\mathcal{D}_{1}\cdot\mathcal{D}_{2})(d,u)&=&  \int_{0}^{\infty} e^{-\frac{t}{d}}h(t,u) d\,t=\int_{0}^{\infty} e^{-\frac{t}{d}}\tilde{h}(t,u) d\,t+\mathcal{O}\Big(d^ke^{-\frac{\pi^2}{d}}\Big) \\
&=&d^k\int_{0}^{\infty} e^{-\frac{t}{d}}\tilde{h}^{(k)}(t,u) d\,t+\mathcal{O}\Big(d^ke^{-\frac{\pi^2}{d}}\Big)\\
&=&d^k\mathcal{D}(d,u)+\mathcal{O}\Big(d^ke^{-\frac{\pi^2}{d}}\Big)
\end{eqnarray*}
where $\mathcal{D}(d,u)$ has property \textmd{G}.
\end{enumerate}
 For $0\leq u\leq 1,$ the proof of  (3), (4), (5) and (6)  is exactly the one of  \cite{SV}. This proof is valid for  $u_0<u\leq 1.$ 

%\section{Approximate solution $ \mathcal{A}(d,u)$ of the functionalequation (\ref{eq1})}

Now we have a formal solution of the equation (\ref{eq1}) and an asymptotic estimate for its coefficients. With the results on the functions in the beginning of this section we have enough information to be able to give a precise function  which satisfies the equation (\ref{eq1}) with an exponentially small error  as $d\to 0$.

%In this section we use the formal solution of  (2.4), the estimates of 
%section 5 for its coefficients and the preliminary results of section
%6 to construct function  $ \mathcal{A}(d,u)$ that almost satisfy (2.4)
%except an error which is exponentially small as $d\to 0$.
In Theorem \ref{the1} we found that (\ref{eq1}) has a uniquely determined formal
power series solution 
\begin{eqnarray}
A(d,u)&=&U(d,u)+Q(d,u)G(d,u)
\labell{re49}\end{eqnarray}
 where $U, Q$ are defined in (\ref{re8}),  (\ref{re10})  and
$G$ is given by (\ref{re36}).
This suggest that we put
\begin{eqnarray}
\mathcal{G}(d,u)&=&\Big(\frac{\alpha}{d^2}+(\beta+\frac{\alpha}{3})\Big)\dfrac{u\, \mathcal{H}_{0}(d,u)}{\tau_2(u)}-\big(\beta d+\frac{\alpha}{d}\big)\mathcal{H}_{2}(d,u)\nonumber\\
&+&\delta d \mathcal{H}_{1}(d,u)+d^{^-2}\mathcal{R}(d,u)\nonumber\\
\mathcal{A}(d,u) &=& U(d,u)+ Q(d,u)\mathcal{G}(d,u)
\labell{re52}\end{eqnarray}
for $d>0, u_0<u\leq1$, 
where $\mathcal{H}_{i}(d,u), i=0..2$ are defined in (\ref{re44}) and
$\mathcal{{R}}(d,u)$ is the function corresponding to
$R_{n},\  n=8.10..$ according to lemma \ref{lem1}.  Using (1) of proposition \ref{propo2} and (2) of lemma \ref{lem1}, we obtain
\begin{equation}
\mathcal{G}(d,u)\sim G(d,u)\quad \textrm {as} \ d\searrow 0\  \textrm
{for every}\,  u_0<u\leq1.
\labell{re53}\end{equation}
Consequently 
\begin{eqnarray}
\mathcal{A}(d,u)&\sim &A(d,u)\quad \textrm {as} \ d\searrow 0\  \textrm
{for every}\,  u_0<u\leq 1.
\labell{re54}\end{eqnarray}
\begin{theo} The function $\mathcal{G}(d,u)$ satisfies (\ref{re16}) except for an exponentially small error. More precisely
\begin{equation*}
\Bigg|V_{1}\cdot\mathcal{S}\Big(Q_{1}\mathcal{S}(\mathcal{G})\Big)-W_{2}\cdot\mathcal{S}^2(\mathcal{G})-W_{3}\mathcal{S}\mathcal{C}(\mathcal{G})-Q_{1}W_{1}\mathcal{G}-Q_{1}f_{2}\Big(d,u,\mathcal{G}\Big)\Bigg|\leq Kd^3(1-u^2) e^{-\frac{\pi^{2}}{d}},
\end{equation*}
uniformly for   $ ( u_0<u<1, 0<d<d_0,)$,
where $K$ is a constant independent of $d$ and $u$.
\labell{the10}\end{theo}

In the proof of this Theorem the functions $\mathcal{D}_i(d,u), i=1.2..$ have property \textmd{G}.
\noindent\textbf{Proof.}
 We set 
\begin{equation}
\mathcal{F}(d,u)=V_{1}\cdot\mathcal{S}\Big(Q_{1}\mathcal{S}(\mathcal{G})\Big)-W_{2}\cdot\mathcal{S}^2(\mathcal{G})-W_{3}\mathcal{S}\mathcal{C}(\mathcal{G})-Q_{1}W_{1}\mathcal{G}-Q_{1}f_{2}\Big(d,u,\mathcal{G}\Big) .
\labell{re55}\end{equation}
Using (2), (4), (5) of  Lemma \ref{lem6.1}, and (\ref{re52}), we obtain 
%According to  (\ref{re46}),  (\ref{re47}) and (\ref{re52}) 
\begin{eqnarray*}
\mathcal{S}(\mathcal{G})=-\Big(\alpha_{1}\sinh(d)\dfrac{1}{\tau_2}-  (\frac{\alpha_1}{2}\sinh(d)-\frac{\beta_1}{2})\Big) \mathcal{H}_{0}+(1-u^2)d^8\mathcal{D}_1(d,u)+d^{-2}\mathcal{S}(\mathcal{R})
\labell{re56}\end{eqnarray*}
where $\mathcal{D}_1(d,u)$ has  property \textmd{G} and

\begin{eqnarray*}
\alpha_1 :&=&\frac{\alpha}{d^2}+(\beta+\frac{\alpha}{3})\\
\beta_1 : &=&\frac{\alpha}{d}+ \beta d
\end{eqnarray*}
This with lemma \ref{lem6.1} imply 
\begin{eqnarray*}
V_1\mathcal{S}\Big(Q_{1}\mathcal{S}(\mathcal{G})\Big)&=& d^5\mathcal{D}_2(d,u)\mathcal{H}_{0}+(1-u^2)d^{11}\mathcal{D}_3(d,u)
+d^{-2}V_1\mathcal{S}(Q_1)\mathcal{C}\mathcal{S}(\mathcal{R})\\
&+&d^{-2}V_1\mathcal{C}(Q_1)\mathcal{S}^{2}(\mathcal{R}),
\end{eqnarray*}
where $ \mathcal{D}_2(d,u),  \mathcal{D}_3(d,u) $ have property \textmd{G}.
%$-\Big(\alpha_{1}\sinh(d)\dfrac{1}{\tau_2}-  (\frac{\alpha_1}{2}\sinh(d)-\frac{\beta_1}{2})\Big)
Because 
\begin{eqnarray*}
\mathcal{C}(Q_{1})&=&(1-u^2)d^2\big(1+\mathcal{O}(d^2)\big),\\
\mathcal{S}(Q_{1})&=&(1-u^2)d^3\big(u+\mathcal{O}(d^2)\big),\\
V_1 &=& 1+\mathcal{O}(d^2)
\end{eqnarray*}
it is sufficient to apply  Theorem \ref{th9}, (4) of Lemma \ref{lem11} for $\mathcal{R}$ and (1) of Lemma \ref{lem11} for  $\mathcal{H}_{0}$  we obtain 
\begin{eqnarray*}
V_1\mathcal{S}\Big(Q_{1}\mathcal{S}(\mathcal{G})\Big)&=& d^5\mathcal{D}_2(d,u)\mathcal{H}_{0}+(1-u^2)d^{2}\mathcal{D}_4(d,u)+\mathcal{O}\Big((1-u^2)d^3 e^{-\frac{\pi^2}{d}}\Big)
\end{eqnarray*}
where $  \mathcal{D}_4(d,u) $ has property \textmd{G}.  With (1) of lemma  \ref{lem11} this implies
\begin{eqnarray}
V_1\mathcal{S}\Big(Q_{1}\mathcal{S}(\mathcal{G})\Big)&=& (1-u^2)d^{2}\mathcal{D}_5(d,u)+\mathcal{O}\Big((1-u^2)d^3 e^{-\frac{\pi^2}{d}}\Big)
\labell{re57}\end{eqnarray}

Using the same method for the terms  $W_{2}\cdot\mathcal{S}^2(\mathcal{G})$, $W_{3}\mathcal{S}\mathcal{C}(\mathcal{G})$ and $Q_{1}W_{1}\mathcal{G}$, we obtain 
\begin{eqnarray}
W_{2}\cdot\mathcal{S}^2(\mathcal{G})&=&(1-u^2)d^{10}\mathcal{D}_6(d,u)+\mathcal{O}\Big((1-u^2)d^{10} e^{-\frac{\pi^2}{d}}\Big),\nonumber\\
W_{3}\mathcal{S}\mathcal{C}(\mathcal{G})&=&(1-u^2)d^{9}\mathcal{D}_7(d,u)+\mathcal{O}\Big((1-u^2)d^{10} e^{-\frac{\pi^2}{d}}\Big),\nonumber\\
Q_{1}W_{1}\mathcal{G}&=&(1-u^2)d^{8}\mathcal{D}_8(d,u)+\mathcal{O}\Big((1-u^2)d^{10} e^{-\frac{\pi^2}{d}}\Big).
\labell{re57.a}\end{eqnarray}
To study the term $Q_{1}f_{2}\Big(d,u,\mathcal{G}\Big) $, we first  treat  $\mathcal{G}^2$ and $\mathcal{G}^3$ .
\begin{eqnarray}
\mathcal{G}^2(d,u)&=&\alpha_1^2\dfrac{u^2}{\tau_2^2}\mathcal{H}_0^2+\beta_1^2\mathcal{H}_2^2+\delta^2 d^2\mathcal{H}_1^2+d^{-4}\mathcal{R}^2-2\alpha_1 \beta_1\dfrac{u}{\tau_2}\mathcal{H}_0 \mathcal{H}_2
+2\alpha_1 \delta\dfrac{u}{\tau_2}\mathcal{H}_0 \mathcal{H}_1 \nonumber\\
&+& 2d^{-2}\alpha_1 \dfrac{u}{\tau_2}\mathcal{H}_0 \mathcal{R}-2\beta_1 \delta d\mathcal{H}_1 \mathcal{H}_2-2d^{-2} \beta_1\mathcal{R} \mathcal{H}_2+2d^{-1}\delta \mathcal{R} \mathcal{H}_1
\end{eqnarray}
Using (1), (2) and (4) of the lemme \ref{lem11} we obtain 
%$$Q_{1}f_{1}\Big(d,u,\mathcal{G}\Big) $$
\begin{eqnarray}
\mathcal{G}(d,u)^2=d^{2}\mathcal{D}_9(d,u)+\mathcal{O}\Big(d^2 e^{-\frac{\pi^2}{d}}\Big)
\labell{re58}\end{eqnarray}
With (2) of lemma \ref{lem11} this implies
\begin{eqnarray}
\mathcal{G}(d,u)^3=\mathcal{G}(d,u)\mathcal{G}(d,u)^{2}=d^{8}\mathcal{D}_{10}(d,u)+\mathcal{O}\Big(d^{8} e^{-\frac{\pi^2}{d}}\Big).
\labell{re58.a}\end{eqnarray}
We can rewrite 
\begin{eqnarray}
f_{2}\big(d,u,\mathcal{G}\big)=y_0(d,u)+y_1(d,u)\mathcal{G}^2 f_{1,1}\big(d,u,\mathcal{G}^2)+y_2(d,u)d^2\mathcal{G}^3 f_{1,2}\big(d,u,\mathcal{G}^2\big)
\labell{58.b}\end{eqnarray}
where $f_{1,1},  f_{1,2}$ and  $ y_i,i=0, 1, 2$ are analytic . 
\begin{lem} 
For $i=1,2$, 
\begin{eqnarray}
f_{1,i}(d,u,\mathcal{G}^2)=\int_{0}^{\infty}e^{-\frac{t}{d}}f_i(t,u)d\,t+\mathcal{O}\Big(d^2 e^{-\frac{\pi^2}{d}}\Big)\quad \text{for} (d>0, u_0<u \leq 1)
\end{eqnarray}
where $f_i(.,u)$ is analytic in $|t|<\pi^2$  and continuous in $[0, \pi^2[ $ and $[\pi^2, \infty[ $
\end{lem} 
\noindent {\bf Proof.}
Using (\ref{re58}), we can rewrite  
$$\mathcal{G}^2(d,u)=\int_{0}^{\infty}e^{-\frac{t}{d}}g(t,u)d\,t+\mathcal{O}\Big(d^2 e^{-\frac{\pi^2}{d}}\Big)$$
where $g(.,u)$ is analytic in $|t|<\pi^2$, it is  continuous in $[0, \pi^2[ $ and $[\pi^2, \infty[ $, twice differentiable in $|t|<\pi^2$. We obtain 
$$f_{1,i}(d,u,\mathcal{G}^2)=\sum_{n=1}^{\infty}f_{n,i}(d,u)\mathcal{G}^{2n}$$
where $$f_{n,i}(d,u)=\int_{0}^{\infty}e^{-\frac{t}{d}}\varphi_{n,i}(t,u)d\,t$$
 and $\varphi_{n,i}(t,u)$ are entire functions. 

Using the proof of theorem 5.1 from \cite{CRSS}, we find  
\begin{eqnarray}
f_{1,i}(d,u,\mathcal{G}^2)=\int_{0}^{\infty}e^{-\frac{t}{d}}f_i(t,u)d\,t+\mathcal{O}\Big(d^2 e^{-\frac{\pi^2}{d}}\Big)
\end{eqnarray}
where the  series $$\sum_{n=1}^{^\infty}(\varphi_{n,i}*g^{*n})(t,u)\Big),\ \ \  g^{*n}=g*...*g,  n \text{ times}$$ 
  is uniformly convergent to a function $f_i(t,u)$ analytic in $|t|<\pi^2$  and continuous in $[0, \pi^2[ $ and $[\pi^2, \infty[ $,  satisfies also 
 $$|f_i(t,u)|\leq K \exp(Kt)\ \ \ \text{for} \ \ t\geq 0, u_0<u\leq 1.$$ Then 
$$
f_{1,i}(d,u,\mathcal{G}^2)=\mathcal{Y}_i(d,u)+\mathcal{O}\Big(d^2 e^{-\frac{\pi^2}{d}}\Big)\quad \text{for} (d>0, u_0<u \leq 1),
$$
where $\mathcal{Y}_i(d,u), i=1,2$ have property \textmd{G}.

%$$f_i(t,u)=\sum_{n=1}^{^\infty}(\varphi_{n,i}*g^{*n})(t,u)$$
This lemma  with (\ref{58.b}) and (2) of lemma \ref{lem11} imply
\begin{eqnarray}
Q_{1}f_{2}\Big(d,u,\mathcal{G}\Big) =d^{6}(1-u^2)\mathcal{D}_{10}(d,u)+\mathcal{O}\Big(d^{6}(1-u^2) e^{-\frac{\pi^2}{d}}\Big)
\labell{re58.c}\end{eqnarray}
Combining  (\ref{re57}),  (\ref{re57.a}) and (\ref{re58.c}), we find 

%Observe that $\mathcal{F}_{i}(d,.)$ have no singularity at $u=0$, whereas 
%${\mathcal T}_j^\pm(d,u)$ have a pole there.
%It can be readily calculated that $\mathcal{F}_{i}(d,u)\in d^{2}\mathcal{G}$ for $i=1..3$  and \,$\mathcal{F}_{5}(d,u)\in d\mathcal{G}$.
%Using Theorem (\ref{th9}), part 2.\ of (\ref{lem1}) and Lemma \ref{lem11}, we obtain
% that 
\begin{eqnarray*}
\mathcal{F}(d,u)&=&d^2(1-u^2)\mathcal{D}(d,u)+\mathcal{K}(d,u)
\end{eqnarray*}
where $\mathcal{D}(d,u)$ has property \textmd{G} and
\begin{equation}
\big|\mathcal{K}(d,u)\big|\leq Kd^3(1-u^2)de^{\frac{-\pi^{2}}{d}}\quad
\textrm{for}\quad (0<d<d_0, \,u_0<u<1).
\labell{re59} \end{equation}
%Using (\ref{lem1}), we obtain  
%\begin{equation*}
%\mathcal{F}(d,u)=\mathcal{D}(d,u)+\mathcal{K}(d,u)
%\end{equation*}
%where $\mathcal{D}(d,u)\in \mathcal{G}$. 
Hence
% Now unfortunately
%$\mathcal{H}^{2}$ has not property \textbf{G}, but we have the
%some whate lemma (6.1) for it. Altogether we have proved that 
\begin{equation*}
\mathcal{F}(d,u)=(1-u^2)d^2\int_0^{\infty}e^{-\frac{t}{d}}q(t,u)dt+\mathcal{K}(d,u)
\end{equation*}
where $q(t,u)$ is analytic in $|t|<\pi^{2}$, it is
continuous on $[0,\pi^{2}[$ and $]\pi^{2},\infty[$, 
has a limit as $t\to \pi^{2}$ for every $u_0<u<1$ and
satisfies 
\begin{equation}
|q(t,u)|\leq K\,e^{Kt},
\labell{re60}\end{equation}
with  a constant  $K$ independent of $u$.  If
$q(t,u)=\sum_{n=0}^{\infty}q_{n}(u)t^{n}$ is the power series
of $q(t,u)$ near $t=0$, Watson's lemma  with (\ref{re59}) imply 
\begin{equation*}
\mathcal{F}(d,u)\sim
\sum_{n=0}^{\infty}n!(1-u^2) q_{n}(u)d^{n+3} \quad \textrm {as} \ d\searrow 0\  \textrm
{for every}\,\   u_0<u<1\ \ .
\end{equation*}
On the other hand because of its definition 
\begin{eqnarray*}
\mathcal{F}(d,u)&\sim &V_{1}\cdot\mathcal{S}\Big(Q_{1}\mathcal{S}(\mathcal{G})\Big)-W_{2}\cdot\mathcal{S}^2(\mathcal{G})-W_{3}\mathcal{S}\mathcal{C}(\mathcal{G})\\
&-&Q_{1}W_{1}\mathcal{G}-Q_{1}f_{2}\Big(d,u,\mathcal{G}\Big)= 0+0d+\cdot\cdot\cdot\ \ \ \ \ \ ,
\end{eqnarray*}
%\begin{eqnarray*}
%\mathcal{F}(d,u)&\sim&\frac{1}{2}B(d,T^{+})-\frac{1}{2}B(d,T^{-})-\varepsilon^{2}+B(d,u)^{2}\\
%&\sim & 0+0d+\cdot\cdot\cdot
%\end{eqnarray*}
since the formal series $G$ satisfies (\ref{re16}). This means that all
  $q_{n}\equiv 0$. Thus we obtain for $u_0<u<1$ with  (\ref {re60}) 
\begin{eqnarray*}
(1-u^2)d^2\int_0^{\infty}e^{-\frac{t}{d}}|q(t,u)|dt &\leq & (1-u^2)d^2\int_{\pi^{2}}^{\infty}e^{-\frac{t}{d}}K e^{Kt}dt\\
&\leq & K(1-u^2)d^3e^{-\frac{\pi^{2}}{d}}    \quad (0<d<d_0)
\end{eqnarray*}
and thus 
\begin{equation*}
\big|\mathcal{F}(d,u)\big|\leq K(1-u^2)d^3e^{-\frac{\pi^{2}}{d}}    \quad (0<d<d_0)\ .
\end{equation*}
We have proved that 
\begin{eqnarray*}
\Bigg|V_{1}\cdot\mathcal{S}\Big(Q_{1}\mathcal{S}(\mathcal{G})\Big)&-&W_{2}\cdot\mathcal{S}^2(\mathcal{G})-W_{3}\mathcal{S}\mathcal{C}(\mathcal{G})
-Q_{1}W_{1}\mathcal{G}-Q_{1}f_{2}\Big(d,u,\mathcal{G}\Big)\Bigg|\\ & \leq & K(1-u^2)d^3e^{-\frac{\pi^{2}}{d}},
\end{eqnarray*}
for $0<d<d_0, \,u_0< u< 1$, i.e.\ ${\mathcal Q}(d,u)$ is a quasi-solution of
(\ref{re16}) on this interval. This implies that the function $\mathcal{A}(d,u)$ defined in (\ref{re52}) is a quasi-solution of
(\ref{eq1}), more precisely 
\begin{eqnarray}
\Bigg|\sqrt{\frac{1-(T^{+})^{2}}{1-u^{2}}}\mathcal{A}(d,T^{+})&+&\sqrt{\frac{1-(T^{-})^{2}}{1-u^{2}}}\mathcal{A}(d,T^{-})-2\mathcal{A}(d,u)-f\big(\varepsilon,\mathcal{A}(d,u)\big)\Bigg|\nonumber\\
&\leq& Kde^{-\frac{\pi^{2}}{d}},\labell{quasi}
\end{eqnarray}
%We need in the next section that it is also a quasi-solution on intervals
%$u_0\leq u\leq 1$, where $-1< u_0<0$. This will be shown in the next version of 
%this manuscript.

\section{Distance Between Points of Manifolds}
Clearly, if $q_{\varepsilon} (t)$ is  an exact  solution of the difference equation (\ref{diff}), then $(q_{\varepsilon}(t),p_{\varepsilon}(t))$, where $p_{\varepsilon}(t)=\frac{1}{\varepsilon}\big(q_{\varepsilon}(t)-q_{\varepsilon}(t-\varepsilon)\big)$, is an exact solution of the the system (\ref{int1}).
In the introduction, we have mentioned that the stable manifold $W_{s}^{-}$ of this system at $A=(0,0)$ is parametrized by $t\to (q^-_{\varepsilon}(t),p^-_{\varepsilon}(t))$ and the unstable manifold $W_{u}^{+}$ of (\ref{int1})  at $B=(2\pi,0)$ is parametrized by $t\to (q^{+}_{\varepsilon}(t),p^{+}_{\varepsilon}(t))$,  where $ (q^{-}_{\varepsilon}(t),p^{-}_{\varepsilon}(t))$ is an exact solution of   (\ref{int1}) and $q^{+}_{\varepsilon}(t)=2\pi -q^-_{\varepsilon}(-t), p^{+}_{\varepsilon}(t)=p^-_{\varepsilon}(-t+\varepsilon)$.

In the previous section, we have constructed a quasi-solution $\mathcal{A}(d,u)$  
for equation (\ref{eq1}), i.e.\ it satisfies  this equation with an exponentially small error. %then $\mathcal{A}^{-}(d,u)=2\pi-\mathcal{A}(d,-u)$ is a quasi-solution colse to the exact solution $q^{-}_{\varepsilon}(t)$. 
We denote by $\tilde{W}_{s}$,
$\tilde{W}_{u}$ the manifolds close to $W^{-}_{s}$  respectively  $W^{+}_{u}$  parametrized by  $t\mapsto
\big(\xi_{-}(t),\varphi_{-}(t)\big)$ respectively  $t\mapsto
\big(\xi_{+}(t),\varphi_{+}(t)\big)$, where
$\xi_{-}(t)=\sqrt{1-u(t)^{^2}}\mathcal{A}\big(d,u(t)\big)+q_{0d}(t)$ and $\xi_{+}(t)=2\pi-\xi_{-}(-t)$, 
$\varphi_{-}(t)=\dfrac{1}{\varepsilon}\big(\xi_{-}(t)-\xi_{-}(t-\varepsilon) \big) $,  $\varphi_{+}(t)=\dfrac{1}{\varepsilon}\big(\xi_{+}(t)-\xi_{+}(t-\varepsilon) \big) $.
%$\varphi_{-}(t)=\varphi_{+}(-t)$.
Here and in the sequel, we often omit to indicate the dependence with respect to 
$\varepsilon$ for the sake of simplicity of notation.

We will first show that the
vertical distance between some point $(q_{1},p_{1})$ of the stable manifold $W_{s}^{-}$ and the
manifold $\tilde{W}_{s}$ is exponentially small. For this purpose, we consider the sequence $Z_{n}=(q_{n},p_{n})$ on the stable manifold  $W_{s}^{-}$, defined by  
\begin{eqnarray}
Z_{n+1}=\binom{q_{n+1}}{p_{n+1}}=\phi\binom{q_{n}}{p_{n}}=\binom{q_{n}+\varepsilon p_{n+1}}{p_{n}+\varepsilon \sin(q_{n})},\ \mbox{for}\ n=1,2,3,....
\end{eqnarray}
There is a sequence $t_{n}$ such that 
\begin{eqnarray}
\xi_{-}(t_{n})=q_{n}=\sqrt{1-u(t_{n})^{^2}}\mathcal{A}\big(d,u(t_{n})\big)+q_{0d}(t_{n}).
\labell{re60.0}\end{eqnarray}

The vertical projection of the point $(q_{n},p_{n})$ on the manifold  $\tilde{W}_{s}$ is the point $(q_{n},b_{n})$, where
 $b_{n}=\varphi_{-}(t_{n})=g(q_{n}):=\varphi_{-}\big(\xi_{-}^{-1}(q_{n})\big)$ for 
$n=1,2,3,...$ .
We denote by $\triangle_{\varepsilon}(n)$ the vertical distance between the point  $(q_{n},p_{n})$ and the manifold  $\tilde{W}_{s}$. Then 
\begin{equation}
\triangle_{\varepsilon}(n)=p_{n}-b_{n}.
\labell{re61}\end{equation}
For the $\phi$-image of the point $(q_n,b_n)$ on $\tilde W_s$, we find
%As the point $(q_{n},b_{n})$ is on the manifold $\tilde{W}_{s}$ and $\mathcal{A}(d,u)$ satisfies   equation (\ref{eq1}) except for an exponentially small error,
%the point $(q_{n},b_{n})$ satisfies 
\begin{eqnarray}
\binom{\tilde{q}_{n+1}}{\tilde{b}_{n+1}}&=&\phi\binom{q_{n}}{b_{n}}=\binom{q_{n}+\varepsilon \tilde{b}_{n+1}}{b_{n}+\varepsilon \sin(q_{n})},\nonumber\\
&=&\binom{\xi_{-}(t_{n})+\varepsilon \tilde{b}_{n+1}}{\varepsilon^{-1}\big(\xi_{-}(t_{n})-\xi_{-}(t_{n}-\varepsilon) \big) +\varepsilon \sin\big(\xi_{-}(t_{n})\big)}\nonumber\\
&=&\binom{2\xi_{-}(t_{n})-\xi_{-}(t_{n}-\varepsilon)  +\varepsilon^{2} \sin\big(\xi_{-}(t_{n})\big)}{\varepsilon^{-1}\big(\xi_{-}(t_{n})-\xi_{-}(t_{n}-\varepsilon) \big) +\varepsilon \sin\big(\xi_{-}(t_{n})\big)}
\end{eqnarray}
As the function $\xi_{-}$  satisfies   equation (\ref{diff}) except for an 
exponentially small error because of (\ref{quasi}) and (\ref{re60.0}), we have
\begin{eqnarray}
 \tilde{q}_{n+1}&=&\xi_{-}(t_{n}+\varepsilon)+e_{n+1}(\varepsilon)\nonumber\\
\tilde{b}_{n+1}&=&\varphi_{-}(t_{n}+\varepsilon)+\dfrac{1}{\varepsilon}e_{n+1}(\varepsilon),
\end{eqnarray}
where $|e_{n+1}(\varepsilon)|\leq K \varepsilon\exp(-\dfrac{\pi^{2}}{\varepsilon})$
with some $K$ independent of $\varepsilon$ and $n$.
Since $g\Big(\xi_{-}(t_{n}+\varepsilon)\Big)=\varphi_{-}(t_{n}+\varepsilon)$, we have 
\begin{eqnarray}
\tilde{b}_{n+1} &=&g\Big(\xi_{-}(t_{n}+\varepsilon)\Big)+\dfrac{1}{\varepsilon}e_{n+1}(\varepsilon)\nonumber\\
                     &=&g\Big(\tilde{q}_{n+1}-e_{n+1}(\varepsilon)\Big)+\dfrac{1}{\varepsilon}e_{n+1}(\varepsilon)
\labell{re62}\end{eqnarray}
% In particular 

%{\bf RE62}
%\begin{equation}
%\tilde{b}_{n+1}-g(\tilde{q}_{n+1})=\varphi_{+}(t_{n}+\varepsilon)+\dfrac{1}{\varepsilon}e_{n+1}(\varepsilon)-g\Big(\xi_{+}(t_{n}+\varepsilon)+e_{n+1}(\varepsilon)\Big)
%\labell{re62}\end{equation}
%using Taylor expansion, we find
%\begin{equation}
%\tilde{b}_{n+1}-g(\tilde{q}_{n+1})=\varphi_{+}(t_{n}+\varepsilon)+\dfrac{1}{\varepsilon}e_{n+1}(\varepsilon)-g\Big(\xi_{+}(t_{n}+\varepsilon)\Big)-g'(\theta_{n+1})e_{n+1}(\varepsilon)
%\end{equation}
%where   $\xi_{+}(t_{n}+\varepsilon)<\theta_{n+1}<\xi_{+}(t_{n}+\varepsilon)+e_{n+1}(\varepsilon)$.

%But, $$g\Big(\xi_{+}(t_{n}+\varepsilon)\Big)=\varphi_{+}(t_{n}+\varepsilon),$$  then 

%{\bf re63}
%\begin{equation}
%\tilde{b}_{n+1}-g(\tilde{q}_{n+1})=\dfrac{1}{\varepsilon}\Big(1-\varepsilon g'(\theta_{n+1})\Big)e_{n+1}(\varepsilon)
%\labell{re63}\end{equation}
\noindent On the other hand, using (\ref{re61})  and the definition of $\phi$, we obtain 
\begin{eqnarray}
\binom{\tilde{q}_{n+1}}{\tilde{b}_{n+1}}&=&\phi\binom{q_{n}}{p_{n}-\triangle_{\varepsilon}(n)}=\binom{q_{n}+\varepsilon p_{n}+\varepsilon^{^2} \sin(q_{n})-\varepsilon\triangle_{\varepsilon}(n)}{p_{n}+\varepsilon \sin(q_{n})-\triangle_{\varepsilon}(n)},\nonumber\\
&=&\binom{q_{n+1}-\varepsilon\triangle_{\varepsilon}(n)}{p_{n+1}-\triangle_{\varepsilon}(n)}.
\end{eqnarray}
%this and (\ref{re62}) imply that
%\begin{equation}
%p_{n+1}-\triangle_{\varepsilon}(n)-g\Big(q_{n+1}-\varepsilon\triangle_{\varepsilon}(n)\Big)=e_{n+1}(\varepsilon).
%\end{equation}
With (\ref{re61}) and   (\ref{re62}) this implies
\begin{eqnarray*}
\triangle_{\varepsilon}(n+1)&=&p_{n+1}-b_{n+1}=p_{n+1}-g(q_{n+1})\nonumber\\
&=& \triangle_{\varepsilon}(n)+\tilde{b}_{n+1}-g\Big(\tilde{q}_{n+1}+\varepsilon\triangle_{\varepsilon}(n)\Big)\nonumber\\
&=&\triangle_{\varepsilon}(n)+g\Big(\tilde{q}_{n+1}-e_{n+1}(\varepsilon)\Big)-g\Big(\tilde{q}_{n+1}+\varepsilon\triangle_{\varepsilon}(n)\Big)+\dfrac{1}{\varepsilon}e_{n+1}(\varepsilon).
%&=&\triangle_{\varepsilon}(n)+e_{n+1}(\varepsilon)-\varepsilon g'(\theta_{n+1})\triangle_{\varepsilon}(n)\nonumber\\
%&=& \big(1-\varepsilon g'(\theta_{n+1})\big)\triangle_{\varepsilon}(n)+e_{n+1}(\varepsilon)
\end{eqnarray*}
Using Taylor expansion we obtain 
%\begin{eqnarray*}
%\triangle_{\varepsilon}(n+1)&=&\triangle_{\varepsilon}(n)-\Big(\tilde{q}_{n+1}-e_{n+1}(\varepsilon)\Big)-g\Big(\tilde{q}_{n+1}+\varepsilon\triangle_{\varepsilon}(n)\Big)+\dfrac{1}{\varepsilon}e_{n+1}(\varepsilon).
%\end{eqnarray*}
\begin{eqnarray}
\triangle_{\varepsilon}(n+1) %&=& \triangle_{\varepsilon}(n)+\tilde{b}_{n+1}-g(\tilde{q}_{n+1}\big)-\varepsilon g'(\vartheta_{n+1})\triangle_{\varepsilon}(n)\nonumber\\
%&=&\triangle_{\varepsilon}(n)+e_{n+1}(\varepsilon)-\varepsilon g'(\theta_{n+1})\triangle_{\varepsilon}(n)\nonumber\\
&=& \big(1-\varepsilon g'(\theta_{n+1})\big)\triangle_{\varepsilon}(n)+\dfrac{1}{\varepsilon}\big(1-\varepsilon g'(\theta_{n+1})\big)e_{n+1}(\varepsilon)\labell{delta1}
\end{eqnarray}
where $\tilde{q}_{n+1}-e_{n+1}(\varepsilon)<\theta_{n+1}<\tilde{q}_{n+1}+\varepsilon\triangle_{\varepsilon}(n)$.

Now 
%\begin{equation*}
%g'(\theta_{n+1})=\frac{\varphi'_{-}(\sigma_{n+1})}{\xi'_{-}(\sigma_{n+1})},\quad
%\text{where}\quad \sigma_{n+1}=\xi_{-}^{-1}(\theta_{n+1}),
%\end{equation*}
$g$ is $\varepsilon$-close to the curve $p=-2\sin(q/2)$, hence
\begin{eqnarray*}
g'(\theta_{n+1})&=%&-\tanh\big(\sigma_{n+1}\big)+O(d)=-u(\sigma_{n+1})+O(d),\\
& -\cos\Big(\dfrac{\theta_{n+1}}{2}\Big)+O(d).
\end{eqnarray*}
Thus given  any positive $\mu<\pi$, there is a positive constant $c$
such that for all $q_{1}\leq \mu$, all $n$ and  sufficiently small $d$,
$$1-\varepsilon g'(\theta_{n+1})\geq1+\varepsilon c\ \ .$$
It is now convenient to write (\ref{delta1}) in the form
$$\triangle_{\varepsilon}(n)=\big(1-\varepsilon g'(\theta_{n+1})\big)^{-1}\triangle_{\varepsilon}(n+1) -\frac1\varepsilon e_{n+1}(\varepsilon)$$
As  $\triangle_{\varepsilon}(n)\to 0$ as  $n\to \infty$, this implies that
%\begin{eqnarray*}
%\big|\triangle_{\varepsilon}(n)\big|\leq K e^{-\frac{\pi^{2}}{\varepsilon}} 
%\sum^{\infty}_{k=
%  n}\frac{1}{\prod_{i}^{i=k}\Big(1-\varepsilon g'\big(\theta_{i+1}\big)\Big)},
%\end{eqnarray*}
%where $K$ is a postive constant.
that there is a positive constant $K$ such that 
\begin{eqnarray*}
\big|\triangle_{\varepsilon}(n)\big|&\leq& Ke^{-\frac{\pi^{2}}{\varepsilon}} 
\sum^{\infty}_{k=0
  }(1+\varepsilon c)^{-k}
\end{eqnarray*}
Consequently
\begin{equation}
\triangle_{\varepsilon}(n)=O\bigg( \frac{1}{\varepsilon}\exp\Big(-\frac{\pi^{2}}{\varepsilon}\Big)\bigg).\labell{delta0}
\end{equation}
In particular 
\begin{eqnarray}
\dist_v((q_1,p_1),\tilde W_{s})=O\bigg( \frac{1}{\varepsilon}\exp\Big(-\frac{\pi^{2}}{\varepsilon}\Big)\bigg),\labell{dists}
\end{eqnarray}
where $(q_{1},p_{1})$ is any point on the stable manifold 
$W_{s}^{-}$, provided $q_1\leq\mu<\pi$; here $\dist_v$ denotes the vertical distance.

The estimate (\ref{dists}) can be extended to any $\mu<2\pi$ and a starting
point $(q_1,p_1)$ with $q_1\leq\mu$ in the following way.
The relation (\ref{delta1}) remains valid, only now we just have the existence of 
some constant $c>0$ such that
$1-\varepsilon g'(\theta_{n+1})\geq1-\varepsilon c\ $ for all $n$. 
As system (\ref{int1}) can be regarded as a one-step numerical method for
the system (\ref{int3}) of differential equations and the starting point is
at a distance ${\cal O}(\varepsilon)$ of its solution $(q_0(t),q_0'(t))$,
results on the convergence of one-step methods can be applied and yield that
$q_k=q_0(t_1+(k-1)\varepsilon)+{\cal O}(\varepsilon)$, where $q_0(t_1)=q_1$, uniformly
for  integer $k$, $1\leq k\leq L/\varepsilon$, where $L$ is any positive constant.
We choose $L$ such $q_0(t_1+L)\leq\mu/2<\pi$.
Repeated application of (\ref{delta1}) now gives
$$|\triangle_\varepsilon(n)-\triangle_\varepsilon(n+[L/\varepsilon])|\leq
K e^{-\frac{\pi^2}\varepsilon}\sum_{k=0}^{[L/\varepsilon]}(1-\varepsilon c)^{-k}
\leq \frac{K e^{Lc}}{c\varepsilon}e^{-\frac{\pi^2}\varepsilon}\ \ .$$
To the quantity $\triangle_\varepsilon(n+[L/\varepsilon])$, 
inequality (\ref{delta0}) can be applied, because 
$q_{n+[L/\varepsilon]}\leq\mu/2<\pi$. Thus (\ref{delta0}) and hence also
(\ref{dists}) remain valid also uniformly for $0<q_1\leq\mu$, provided
$\mu<2\pi$.

The analogous reasoning applies to the vertical distance of a point 
$(\tilde q_{1},\tilde p_{1})$ on the unstable manifold $W_u^+$ from the manifold
$\tilde{W}_{u}$ and yields
\begin{eqnarray}
\dist_v((\tilde q_1,\tilde p_1),\tilde W_{u})=O\bigg( \frac{1}{\varepsilon}\exp\Big(-\frac{\pi^{2}}{\varepsilon}\Big)\bigg).\labell{distu}
\end{eqnarray}
Another method to obtain (\ref{distu}) consists in using (\ref{dists}) and symmetry.

Now we will estimate the vertical distance 
between the two manifolds  $\tilde{W}_{s}$ and $\tilde{W}_{u}$.
%We will need to extend both
%quasi-solutions in an adequate way. This is possible because we can show that  for  $-1 \leq u <0$, $\mathcal{A}(d,u)$ satisfies  equation (\ref{eq1}) with an exponentially small error (The demonstration will be sent later).
As the quasi-solution ${\mathcal A}(d,u)$ is defined for $-1<u=:\tanh(t)<1$, we 
can define
\begin{eqnarray}
\mathcal{A}^{+}(d,u)&=&-\mathcal{A}(d,-u),\nonumber\\
\xi_{+}(t)&=&\sqrt{1-u(t)^{2}}\mathcal{A}^{+}\big(d,u(t)\big)+2\pi-q_{0d}(-t)=
2\pi-\xi_-(-t),
\labell{dc01}\end{eqnarray}
\begin{equation}
 D_{\varepsilon}(t)=\xi_{+}(t)-\xi_{-}(t) \quad \text {for} \, -\dfrac{4}{3}<t<\dfrac{4}{3}\ \ .
\labell{dc1}\end{equation}
Using (\ref{dc01}) and the  definition of  $\xi_{-}(t)$ we find 
\begin{eqnarray}
 D_{\varepsilon}(t)&=&-\sqrt{1-u(t)^{^2}}\Big(\mathcal{A}\big(d,u(t)\big)+\mathcal{A}\big(d,-u(t)\big)\Big)-q_{0d}(t)-q_{0d}(-t)+2\pi\nonumber\\
                         &=&-\sqrt{1-u(t)^{^2}}\Big(\mathcal{A}\big(d,u(t)\big)+\mathcal{A}\big(d,-u(t)\big)\Big)
\end{eqnarray}
With (\ref{re13}) and (\ref{re52}) this implies
\begin{eqnarray*}
D_{\varepsilon}(t)&=&-\sqrt{1-u(t)^{^2}}Q(d,u)\Bigg[\alpha_{1}\dfrac{u}{\tau_2}\Big(\mathcal{H}_{0}(d,u)-\mathcal{H}_{0}(d,-u)\Big)
-\nonumber \\
&&\beta_{1}\Big(\mathcal{H}_{2}(d,u)+\mathcal{H}_{2}(d,-u)\Big)+ \delta d\Big(\mathcal{H}_{1}(d,u)+\mathcal{H}_{1}(d,-u)\Big)\Bigg]
\end{eqnarray*}
where  $Q(d,u)$ is defined in  (\ref{re10}), $\alpha_{1},  \beta_{1}$ are defind in (\ref{re56}) and 
%\begin{eqnarray*}
%\mathcal{H}_i(d,u):= \int^{+\infty}_0e^{-\frac{t}{d}}h_i(t,u) d\,t ,  i=1..,3
%\end{eqnarray*}
%with $ h_i(t,u),  i=1..,3$ 
${\mathcal H}_i(d,u)$ are defined in (\ref{re44}) and can be continued analytically
to $-1<u\leq1$ as in (\ref{re44a}).

Using the fact that the functions $u\,h_0(s,u), h_1(s,u), h_2(s,u)$ in (\ref{re44}) are odd,
we can apply the residue theorem and obtain for $-1<u<1$ that
\newcommand{\Res}{\mbox{\,Res}}
\begin{eqnarray*}
\mathcal{H}_{i}(d,u)+\mathcal{H}_i(d,-u)&=&
\sum_{Im(s_{k}(t))<0}2\pi i
\Res\Big(e^{-\frac{s}{d}}h_{i}(s,u),s_{k}(t)\Big)\\
&&-\sum_{Im(s_{k}(t))>0}2\pi i
\Res\Big(e^{-\frac{s}{d}}h_{i}(s,u),s_{k}(t)\Big),\quad i=0,1,2,
\end{eqnarray*}
where $s_{k}(t)= \pi^{2}\pm
\frac{2d\pi\,t}{\varepsilon}\,i+2k\pi^{2}$ \,for $k\ge 0$.
%  $\Gamma$ is
%the sum of two paths $\Gamma_{1}$ and  $\Gamma_{2}$. Here
%$\Gamma_{1}$ consist of two segments of the higher half-plane, one of
%these two segments is parallel to the axis $y=0$ and begins at the  point $(\frac{\pi^{2}}{4},1)$, the other ends
%in the point (0,0), and  the path $\Gamma_{2}$  is the symmetry of
%$\Gamma_{1}$ relatively to the axis $y=0$.
We obtain 
\begin{eqnarray*}
\Res\Big(e^{-\frac{s}{d}}h_{0}(s,u),s_{k}(t)\Big)&=&\dfrac{1}{2}e^{-\dfrac{(k+1)\pi^{2}}{d}}e^{\mp\dfrac{2\pi t i}{\varepsilon}}\nonumber\\
\Res\Big(e^{-\frac{s}{d}}h_{1}(s,u),s_{k}(t)\Big)&=&\pm\dfrac{i}{4\pi}e^{-\dfrac{(k+1)\pi^{2}}{d}}e^{\mp\dfrac{2\pi t i}{\varepsilon}}\nonumber\\
%&\pm &\dfrac{i}{2\pi d}e^{-\dfrac{(k+1)\pi^{2}}{d}}e^{\mp\dfrac{2\pi t i}{\varepsilon}}\nonumber\\
\Res\Big(e^{-\frac{s}{d}}h_{2}(s,u),s_{k}(t)\Big)&=&\pm \dfrac{i}{4\varepsilon d}\big(\pi \varepsilon i \mp 2d t\big)e^{-\dfrac{(k+1)\pi^{2}}{d}}e^{\mp\dfrac{2\pi t i}{\varepsilon}}
\end{eqnarray*}
and hence
\begin{eqnarray}
D_{\varepsilon}(t)&=&\frac{1}{d^2} \phi_{1}(t,\varepsilon)\exp\bigg(-\frac{\pi^{2}}{d}\bigg)+O\bigg(e^{-\frac{\pi^{2}}{d}}\bigg).
\labell{dc2}\end{eqnarray}
where
\begin{eqnarray}
\phi_{1}(t,\varepsilon)&=&2\pi\alpha\bigg[\sinh\Big(\frac{dt}{\varepsilon}\Big)+t/\cosh\Big(\frac{dt}{\varepsilon}\Big)\bigg]\sin\Big(\frac{2\pi t}{\varepsilon}\Big)\nonumber\\
&+&\bigg[\dfrac{ \pi \alpha}{\cosh\Big(\frac{dt}{\varepsilon}\Big)}\bigg]\cos\Big(\frac{2\pi t}{\varepsilon}\Big)
%\phi_{2}(t,\varepsilon)&=&-\bigg[\dfrac{4\beta_{2}(\varepsilon)}{\cosh\Big(\frac{dt}{\varepsilon}\Big)}-\gamma_{2}(\varepsilon)\cosh\Big(\frac{dt}{\varepsilon}\Big)\bigg]\cos\Big(\frac{2\pi t}{\varepsilon}\Big)\nonumber
\end{eqnarray}
As a consequence of (\ref{dc2}) and (\ref{dc01}), we obtain immediately that
\begin{equation}\xi^\pm(0)=\pi+{\cal O}(\varepsilon^{-2}e^{-\pi^2/\varepsilon})\ \ .
\labell{xi0pi}\end{equation}

%\begin{eqnarray}
%d_{\varepsilon}(t)&=&\frac{\alpha\sqrt{1-u(t)^{^2}}}{d^3} \Bigg[5u(t)\cosh\Big(\frac{dt}{\varepsilon}\Big)\cos\Big(\frac{2\pi t}{\varepsilon}+\pi\Big)
%-5 t\cos\Big(\frac{2\pi t}{\varepsilon}\Big)\nonumber \\
%&-&\frac{5\pi}{2}\sin\Big(\frac{2\pi t}{\varepsilon}\Big)\Bigg]\exp\bigg(-\frac{\pi^{2}}{d}\bigg)+O\bigg(\frac{1}{d^{3}}e^{-\frac{\pi^{2}}{d}}\bigg).
%\labell{dc2}\end{eqnarray}
Now, let us take  a point $(\xi_{+}(t), \varphi_{+}(t)) $ on the  
manifold $\tilde{W}_{u}$. We suppose that the point $(\xi_{-}(t_1), \varphi_{-}(t_{1}))$ is its vertical projection on the manifold 
$ \tilde{W}_{s} $. We will evaluate the  vertical distance  between these two points
\begin{eqnarray}
\dist_{v}(t)=\varphi_{+}(t)-\varphi_{-}(t_{1})=\varphi_{+}(t)-g\big(\xi_{+}(t)\big),
\end{eqnarray}
where $g(x)=\varphi_{-}\Big(\xi_{-}^{-1}(x)\Big)$.
Thus by (\ref{dc1})
\begin{eqnarray}
\dist_{v}(t)=\varphi_{+}(t)-g\big(\xi_{-}(t)+D_{\varepsilon}(t)\big).
\end{eqnarray}
Using Taylor expansion, we find 
\begin{eqnarray}
\dist_{v}(t)=\varphi_{+}(t)-g\big(\xi_{-}(t)\big)-D_{\varepsilon}(t)g'\big(\eta(t)\big),
\end{eqnarray}
where $\xi_{-}(t)<\eta(t)<\xi_{-}(t)+D_{\varepsilon}(t)$. Here
\begin{eqnarray}
\eta(t)&=&q_{0d}(t)+O(d) ,\ \ \mbox{hence}\\
g'\big(\eta(t)\big)&=&-\cos\Big(\dfrac{\eta(t)}{2}\Big)+O(d)=-\tanh\Big(\frac{dt}{\varepsilon}\Big)+O(d).                     
\labell{dr0}\end{eqnarray} 
As $g\big(\xi_{-}(t)\big)=\varphi_{-}(t)$, this yields
\begin{eqnarray}
\dist_{v}(t)&=&\varphi_{+}(t)-\varphi_{-}(t)-D_{\varepsilon}(t)g'_{-}\big(\eta(t)\big),\nonumber\\
             &=&\frac{1}{\varepsilon}\big(\xi_{+}(t)-\xi_{+}(t-\varepsilon)\big)-\frac{1}{\varepsilon}\big(\xi_{-}(t)-\xi_{-}(t-\varepsilon)\big)-D_{\varepsilon}(t)g'_{-}\big(\eta(t)\big),\nonumber \\
&=&\frac{1}{\varepsilon}\big(D_{\varepsilon}(t)-D_{\varepsilon}(t-\varepsilon)\big)-D_{\varepsilon}(t)g'\big(\eta(t)\big).
\labell{dr1}\end{eqnarray}
Now formula (\ref{dc2}) applies and we obtain
\begin{eqnarray}
D_{\varepsilon}(t-\varepsilon)&=&D_{\varepsilon}(t)+\dfrac{1}{d}\phi_{2}(t,\varepsilon)\exp\bigg(-\frac{\pi^{2}}{d}\bigg)+O\bigg(e^{-\frac{\pi^{2}}{d}}\bigg)
\end{eqnarray}
where
\begin{eqnarray}
\phi_{2}(t,\varepsilon)&=&-2\pi\alpha\bigg[\dfrac{\cosh\Big(\frac{dt}{\varepsilon}\Big)^2+1}{\cosh\Big(\frac{dt}{\varepsilon}\Big)}-t\dfrac{\,\tanh\Big(\frac{dt}{\varepsilon}\Big)}{\cosh\Big(\frac{dt}{\varepsilon}\Big)}\bigg]\sin\Big(\frac{2\pi t}{\varepsilon}\Big)\nonumber\\
&+& \pi\alpha\bigg[\dfrac{\tanh\Big(\frac{dt}{\varepsilon}\Big)}{\cosh\Big(\frac{dt}{\varepsilon}\Big)}\bigg]\cos\Big(\frac{2\pi t}{\varepsilon}\Big)\ \ .
\end{eqnarray}
With (\ref{dr1}) and  (\ref{dr0}) this implies
\begin{eqnarray}
\dist_{v}(t)=\dfrac{1}{d^{2}}\bigg[-\phi_{2}(t,\varepsilon)+\tanh\Big(\frac{dt}{\varepsilon}\Big)\phi_{1}(t,\varepsilon)\bigg]e^{-\frac{\pi^{2}}{d}}+O\bigg(\frac{1}{d}e^{-\frac{\pi^{2}}{d}}\bigg)\ \ .
\end{eqnarray}
Consequently, for $-\dfrac{4}{3}\leq t\leq \dfrac{4}{3}$
\begin{eqnarray}
\dist_{v}(t)=\dfrac{4\pi\alpha}{\varepsilon^{2}}\cosh(t)\sin\Big(\frac{2\pi t}{\varepsilon}\Big)e^{-\frac{\pi^{2}}{\varepsilon}}+O\bigg(\frac{1}{\varepsilon}e^{-\frac{\pi^{2}}{\varepsilon}}\bigg)\quad \text{as} \,\,\,\varepsilon\searrow 0.
\end{eqnarray}

Combining this result with (\ref{dists}) and (\ref{distu}), we finally obtain
our main result theorem 1.1, because $\xi_\pm(t)=q_0(t)+{\cal O}(\varepsilon^2)$
uniformly with respect to $t$ on any finite interval.
% as (\ref{xi0pi}) and the choice of the parametrisation
%of the unstable manifold imply that points with same projection on the $x-axis$
%have  

%\begin{eqnarray}
%\dist_{v}(t)=\dfrac{\alpha\varphi_0(t,\varepsilon)}{\varepsilon^{2}}e^{-\frac{\pi^{2}}{\varepsilon}}+O\bigg(\frac{1}{\varepsilon}e^{-\frac{\pi^{2}}{\varepsilon}}\bigg)
%\end{eqnarray}
%where 
%\begin{eqnarray}
%\varphi_0(t,\varepsilon)&=&\dfrac{\pi}{2}\dfrac{\cosh(t)+(1-2t)\sinh(t)}{\cosh(t)^{2}}\sin\Big(\frac{2\pi t}{\varepsilon}\Big)\nonumber\\
%&+&\dfrac{\pi^{2}}{4}\dfrac{ \sinh(t)}{\cosh(t)^{2}}\cos\Big(\frac{2\pi t}{\varepsilon}\Big)
%\end{eqnarray}

 %$$\varphi_{+}(t)=\frac{1}{\varepsilon}\big(\xi_{+}(t)-\xi_{+}(t-\varepsilon)\big),$$ then 
%{\bf DR2}
%\begin{equation}
%\varphi_{+}(t+\varepsilon)-\varphi_{+}(t)=\frac{1}{\varepsilon}\big(\xi_{+}(t+\varepsilon)-2\xi_{+}(t)+\xi_{+}(t-\varepsilon)\big).
%\labell{dr2}\end{equation}
%On the other hand, by the construction of the function $f$ defined in (\ref{eq1}) and  theorem \ref{the10} and   the definition of $\xi_{+}(t)$ in (\ref{re60.0}),  $\xi_{+}(t)$ satisfies equation (\ref{int0}) with an exponentially small error, and Similar for the function $\xi_{-}(t)$. Precisely
%$\substack{n=1\\ n\ odd}$

%\newpage


\begin{thebibliography}{02}
\bibitem{CRSS}
{M. Canalis-Durand, J. P. Ramis, R. Sch\"afke, Y. Sibuya},  Gevrey solutions of singularly perturbed differential equations,
{\em J. reine. Angew. Math.} 518 (2000), 95-129.
\bibitem{WE}
{W. Eckhaus}, Asymptotic Analysis of Singular perturbations,
{\em North-Holland, Amsterdam} (1979).
\bibitem{FS}
{A. Fruchard, R. Sch\"afke}, Exponentially small splitting of
separatrices for difference equations with small step size, {\em Journal
of Dynamical and Control Systems} 2 (1996), no. 2, 193-238.
\bibitem{G}
{V. G. Gelfreich}, A proof of the exponentially small transversality of
the separatrices for the standard map.{\em Comm. Math. Phys.} 201 (1999), no.
1, 155-216.
\bibitem{GL}
{V. G. Gelfreich, V. F. Lazutkin, and N. V. Svanidze.}  A refined formula for the separatrix splitting for the standard map.{\em Physica D} 71(2), 82-101 (1994)
\bibitem{HM}
{V. Hakim, K. Mallick}, Exponentially small splitting of separatrices,
matching in the complex plan and Borel summation,{\em Nonlinearity} 6 (1993) 
57-70.
\bibitem{L1}
{V.F. Lazutkin}, Splitting of separatrices for the Chirikov's standard map.  VINITI no. 6372/84, (1984), 
(Russian)
\bibitem{L2}
{V.F. Lazutkin}, Exponential Splitting of
separatrices and an analytical integral for the semistandard
map. {\em Preprint, Université Paris VII}, (1991) cf. MR 94a:58108.  
\bibitem{LS}
{V.F. Lazutkin, I.G. Schachmannski and M.B. Tabanov}, Splitting of
separatrices for standard and semistandard mappings, {\em Physica D,}
40 235-248, (1989).
\bibitem{SV}
{R. Sch\"afke, H. Volkmer}. Asymptotic analysis of the equichordal problem,
{\em J. reine . Angew. Math.} 425 (1992), 9-60.
\bibitem{S}
{ H. Sellama}, On the distance between separatrices for the discretized 
logistic differential equation, submitted. http://hal.archives-ouvertes.fr/docs/00/28/75/86/PDF/article-logistic.pdf.
\bibitem{SY}
{Suris, Yuri B.} On  the complex separatrices of some standard-like maps. {\em Nonlinearity} 7 (1994), no. 4, 1225-1236.
\end{thebibliography}
\end{document}